\documentclass{article}
\usepackage[ansinew]{inputenc}
\usepackage[english]{babel}
\usepackage{graphicx}
\usepackage{amsmath,amssymb,amsthm,mathrsfs,amsfonts}
\usepackage{cancel}

\newtheorem {thm}{Theorem}[section]
\newtheorem {defn}{Definition}[section]

\newtheorem {mthm}[thm]{Main Theorem}
\newtheorem {lemma}[thm]{Lemma}

\newtheorem {cor}[thm]{Corollary}

\newcommand{\I}{\boldsymbol{\mathcal{I}}}

\newcommand{\SU}{\boldsymbol{\mathcal{SU}}}
\newcommand{\U}{\boldsymbol{\mathcal{U}}}
\newcommand{\R}{\boldsymbol{\mathcal{R}}}

\newcommand{\RSU}{\boldsymbol{r\mathcal{SU}}}
\newcommand{\RU}{\boldsymbol{r\mathcal{U}}}
\newcommand{\RR}{\boldsymbol{r\mathcal{R}}}

\newcommand{\SP}{{\operatorname{sp}}}
\newcommand{\TR}{{\operatorname{tr}}}
\newcommand{\SPEC}{{\operatorname{spec}}}

\newcommand{\Rem}{\noindent\textbf{Remark: }}
\newcommand{\Rems}{\noindent\textbf{Remarks: 1. }}
\newcommand{\Proof}{\noindent\textbf{Proof: }}

\newcommand{\bbmatrix}[1]{\left[ \begin{array}{cccccccccccccccccc} #1 \end{array} \right]}

\title{Solution of the KdV equation on the line with analytic initial potential}
\author{Andrey Melnikov\\Drexel University, Philadelphia, USA}
\begin{document}
\maketitle

\abstract{We present a theory of Sturm-Liouville non-symmetric vessels, realizing an inverse scattering theory for the Sturm-Liouville operator
with analytic potentials on the line. This construction is equivalent to the construction of a matrix spectral measure for the Sturm-Liouville operator, defined
with an analytic potential on the line.
Evolving such vessels we generate KdV vessels, realizing solutions of the KdV equation.
As a consequence, we prove the following theorem:

Suppose that $q(x)$ is an analytic function on $\mathbb R$. There exists a KdV vessel, which exists on $\Omega\subseteq \mathbb R^2$.
For each $x\in\mathbb R$ there exists $T_x>0$ such that $\{x\}\times[-T_x,T_x]\in\Omega$. The potential $q(x)$ is realized by the vessel for $t=0$. 

Since we also show that if $q(x,t)$ is a solution of the KdV equation on $\mathbb R\times[0,t_0)$, then there exists a
vessel, realizing it, the theory of vessels becomes a universal
tool to study this problem. Finally, we notice that the idea of the proof applies to a similar existence of a solution for evolutionary NLS and Boussinesq equations, since
both of these equations possess vessel constructions.
}

\tableofcontents

\section{Introduction}
The Korteweg-de Vries (KdV) is the following nonlinear evolutionary Partial Differential Equation (PDE) for a
function of two real variables $q(x,t)$:
\begin{equation} \label{eq:KdV}
q_t = -\dfrac{3}{2} q q_x + \dfrac{1}{4} q_{xxx},
\end{equation}
where $q_t, q_x$ denote the partial derivatives. The equation is named after Diederik Korteweg and Gustav de Vries who studied it in \cite{bib:KdV}.
Usually, one considers the \textit{initial value problem}, which is defined as follows: find a solution $q(x,t)$ of \eqref{eq:KdV}, which additionally satisfies:
$q(x,0) = q(x) (x\in\mathbb R)$, for a given function $q(x)$, defined on $\mathbb R$.
A standard technique to solve KdV involves a more elementary equation, called Sturm Liouville (SL) differential equation:
\begin{equation}\label{eq:SL}
-\frac{d^2}{dx^2} y(x) + q(x) y(x) = \lambda y(x),
\end{equation}
where $\lambda\in\mathbb C$ is called the \textit{spectral parameter} and the coefficient $q(x)$ is called the \textit{potential}. In order to solve
\eqref{eq:KdV} using \eqref{eq:SL} one transforms \cite{bib:GGKM} the potential $q(x)$ appearing in \eqref{eq:SL} to 
its ``scattering data". Then one evolves with $t$ the scattering data using some simple differential equations.
Finally, transforming back the evolved scattering data we obtain a new potential $q(x,t)$ of two variables, which solves
\eqref{eq:KdV} and satisfies $q(x,0)=q(x)$. In other words, in this manner we
solve the initial value problem for the equation \eqref{eq:KdV}. So, in order to solve \eqref{eq:KdV} one
has to find a "scattering data" for the given potential.

On the half line the question of characterizing of scattering data (or more precisely spectral measure 
$d\rho(\mu)$ on $\mathbb R$) for
a given potential was completely solved for a continuously differentiable potential by Gelfand-Levitan theory
\cite{bib:GL}, but not always it can be used to solve the KdV equation \eqref{eq:KdV}. For this situation in case 
$var[d\rho]<\infty$ a solution of KdV in the first quadrant ($x,t\geq 0$) is presented in \cite{bib:UnboundedVessels}.

Although there is a good scattering theory of the SL equation \eqref{eq:SL} on the line
\cite[Chapter 2, p. 128]{bib:MarchenkoSL}, \cite[Chapter 2]{bib:LevZar}, even with arbitrary singularities 
\cite{bib:DeiftZho}, the solutions of KdV, corresponding to them are not developed.
In fact, the classes of initial potentials, for which solutions of \eqref{eq:KdV} were presented using inverse scattering are as follows:
\begin{enumerate}
	\item Soliton solutions correspond to $d\rho(\mu)$ to be a finite sum of point mass measures (discrete measure) \cite{bib:Crum},
	\item Krichever solution \cite{bib:Krich77}, where $d\rho(\mu)=f(\mu)d\mu$ for $\mu\in\Gamma$ 
		for some algebraic curve $\Gamma$,
	\item Fadeyev inverse scattering theory \cite{bib:FaddeyevII} where the $d\rho(\mu)$ is supported on 
		the positive real line and has a finite number of point-mass measures on the negative real line,
	\item Periodic potentials \cite{bib:MagWin} 
		correspond to discrete spectral measures with accumulation point at infinity,
	\item Quasi-periodic potentials \cite{bib:BJ45, bib:DSinai}.
\end{enumerate}
While analyzing the KdV equation \eqref{eq:KdV} with an analytic initial potential, one can apply Taylour series techniques to try to solve the equation.
It turns out that the corresponding combinatorial problem is extremely dificult. In fact a recent result of M. Goldstein and D. Damanik \cite{bib:GoldDam} proving 
existence of a global solution of the KdV equation \eqref{eq:KdV} with a quasi periodic potential involves an extremely sophisticated combinatorics of the powers of exponents,
corresponding to all ``harmonics''. Still, the general problem of constructing a solution of the KdV equation for a given analytic potential has yet to be solved.

Using theory of vessels, we show that it is indeed a rare case that a solution of the KdV equation would exists on $\mathbb R\times[0,t_0)$ for some $t_0$.
And the reason for this is that there is an operator of the form $\mathbb X(x,t)=I+T(x,t)$, $I$ - identity, $T(x,t)$ - trace class, 
which is usually invertible for $t=0$ for all $x$, but fails to
be invertible uniformly for all $x$ for whatsoever $t>0$ is. Theoretically it explains why there is no a theory on the existence of the local solution of the KdV equation for arbitrary
analytic initial potential and only special cases ares solvable: for example, in the Faddeyev case the inverse of $\mathbb X(x,t)$ is uniformly bounded, so the continuous
perturbation of the inverse will exist on $[0,t_0)$ (see Remarks). In fact, the existence of the solution $q(x,t)$ of \eqref{eq:KdV} on $\mathbb R\times[0,t_0)$
implies that there exists a vessel on the same set (see Theorem \ref{thm:XinvForq}). This actually means that the theory of vessels is a universal tool to study solutions
of the KdV equation \eqref{eq:KdV}.

We present now the Main Theorem.

\noindent\textbf{Main Theorem \ref{mthm}}
\textit{Suppose that $q(x)$ is an analytic function on $\mathbb R$. There exists a KdV vessel, which exists on $\Omega\subseteq\mathbb R^2$. For each $x\in\mathbb R$ there exists $T_x>0$ such that $\{x\}\times[-T_x,T_x]\in\Omega$. The potential $q(x)$ is realized by the vessel for $t=0$.}

The idea of the proof for this Theorem is simple and appears in \cite{bib:KdVVessels}. Using simple algebra calculation, it is possible to show that
constructing a collection $\mathfrak{V}_{KdV}$ of bounded operators and spaces ($\mathcal H$ - Hilbert space)
\[ \begin{array}{lllllllll}
\mathfrak{V}_{KdV} = \bbmatrix{C(x) & A_\zeta, \mathbb X(x), A & B(x)&  \sigma_1,\sigma_2,\gamma,\gamma_*(x) \\
\mathbb C^2 & \mathcal{K} & \mathbb C^2 & \Omega}, \\
\hspace{2cm} B, C^*:\mathbb C^2\rightarrow\mathcal H, \quad A,\mathbb X,A_\zeta:\mathcal H\rightarrow\mathcal H, \\
\hspace{2cm} \sigma_1, \sigma_2, \gamma,\gamma_*(x, t) \text{ - $2\times 2$ matrices}, \quad \Omega\subseteq\mathbb R^2\\
\dfrac{\partial}{\partial x} B \sigma_1 = - (A B \sigma_2 + B \gamma), \hspace{1.2cm} \dfrac{\partial}{\partial t} B = i A  \dfrac{\partial}{\partial x} B, \\
\dfrac{\partial}{\partial x} C^* \sigma_1  = - (A_\zeta^* C^* \sigma_2 + C^* \gamma), \quad \quad \dfrac{\partial}{\partial t} C^*  = i A_\zeta^*  \dfrac{\partial}{\partial x} C^*, \\
\dfrac{\partial}{\partial x} \mathbb X =  B \sigma_2 C, \hspace{3cm} \dfrac{\partial}{\partial t} \mathbb X = i A B \sigma_2 C - i B \sigma_2 C A_\zeta + i B \gamma C, \\ 
A \mathbb X + \mathbb X A_\zeta  + B \sigma_1 C  = 0, \\
\gamma_* =  \gamma + \sigma_2 C \mathbb X^{-1} B \sigma_1 - \sigma_1 C \mathbb X^{-1} B \sigma_2.
\end{array} \]
which is called a \textbf{regular KdV vessel} on $\Omega$, where $\mathbb X(x, t)$ is invertible, we obtain a solution of \eqref{eq:KdV} on $\Omega$ as follows. For 
\[ \sigma_1 = \bbmatrix{0&1\\1&0}, \quad \sigma_2=\bbmatrix{1&0\\0&0}, \quad \gamma=\bbmatrix{0&0\\0&i}
\]
the function ($(x_0,t_0)\in\Omega$)
\[ q_{\mathfrak V}(x,t) = -2\dfrac{\partial^2}{\partial x^2} \ln \det(\mathbb X^{-1}(x_0,t_0)\mathbb X(x,t))
\]
is analytic in both variables on $\Omega$ and satisfies \eqref{eq:KdV}. The main contribution of this paper is that one can construct a vessel $\mathfrak V_{KdV}$ such that
for $t=0$ it holds that $q_{\mathfrak V}(x,0) = q(x)$ for arbitrary given analytic function $q(x)$. We use unbounded operators on Krein spaces in this case. Moreover, the operators
are $A=i\mu$ and $A_\zeta=-i\mu-\int_0^\infty d\bar\rho(\mu)\cdot$ for a $2\times 2$ measure $d\bar\rho$ on $[0,\infty)$, 
creating a Krein space $\mathcal K=L^2(d\bar\rho)$ (see Section \ref{sec:Back} for details), used instead of the 
Hilbert space $\mathcal H$ in the original definition. The main ingredients of this construction are node, prevessel and vessel. 
A node is the $(0,0)$ value of the vessel $\mathfrak V_{KdV}$, and prevessel is a globally defined object (for all $x,t\in\mathbb R$), which does not include the matrix $\gamma_*(x,t)$
(thus there is no need to demand existence of the inverse of $\mathbb X(x,t)$).
These notions are fully studied at the text along with their properties. Finally, the tau function
$\tau(x,t) = \det (\mathbb X^{-1}(0,0)\mathbb X(x,t))$ defines the set $\Omega\subseteq\mathbb R^2$ where the vessel $\mathfrak V_{KdV}$ exists and, as a result, where the
the solution of \eqref{eq:KdV} exists as well. We show in Main Theorem \ref{mthm} that $\{x\}\times[-T_x,T_x]\subseteq\Omega$ for some $T_x>0$, depending on $x$.

From the Main Theorem \ref{mthm} it follows that it is difficult to create a general construction of solutions of \eqref{eq:KdV} on a strip
$\mathbb R\times[0,t_0)$. Instead, the focus must be made on the problem of characterization of classes of functions $q(x)$, for which
the matrix $\mathbb X(x,t)$ is invertible. The formula for the inverse, if it exists is also given in terms of the fundamental solutions of \eqref{eq:KdV}.
It seems that the ideas, presented in this text also have a potential to be applied to the locally integrable case: $\int_0^x q(y)dy<\infty$ for all $x\in\mathbb R$.

Similar calculations can be carried over for the vessel parameters 
\[ \sigma_1 = I= \bbmatrix{1&0\\0&1}, \quad \sigma_1 = \dfrac{1}{2}\bbmatrix{1&0\\0&-1}, \quad \gamma = 0 = \bbmatrix{0&0\\0&0},
\]
which create solutions $y=\bbmatrix{0&1}\gamma_*(x,t)\bbmatrix{1\\0}$ of the evolutionary Non Linear Schr\" odinger (NLS) equation
\begin{equation} \label{eq:ENLS}
i y_t + y_{xx} + 2 |y|^2 y = 0.
\end{equation}
The proof of this fact can be found in \cite{bib:ENLS}. As a result, the ideas presented in this work can be used to prove a similar to Main Theorem \ref{mthm} result for
the evolutionary NLS equation.

Defining
\[ \sigma_1 = \bbmatrix{0 & 0 & 1 \\ 0 & 1 & 0 \\ 1 & 0 & 0},
   \sigma_2 = \bbmatrix{1 & 0 & 0 \\ 0 & 0 & 0 \\ 0 & 0 & 0},
   \gamma   = \bbmatrix{0 & 0 & 0 \\ 0 & 0 & 1 \\ 0 & -1 & 0}.
\]
and
\[
\widetilde{\sigma}_1 = \sigma_1, \quad \widetilde\sigma_2=\bbmatrix{0&-i&0\\i&0&0\\0&0&0},\quad \widetilde\gamma=\bbmatrix{0&0&0\\0&0&0\\0&0&i},
\]
it was shown in \cite{bib:BoussVessels} that the collection
\begin{equation*}
\mathfrak{V}_{Bouss} = (C(x,t), A_\zeta, \mathbb X(x,t), A, B(x,t); \sigma_1, 
\sigma_2, \gamma, \gamma_*(x,t),\widetilde{\sigma}_1, \widetilde\sigma_2, \widetilde\gamma ;
\mathcal{H},\mathbb C^3;\Omega),
\end{equation*}
where the operators $C(x,t):\mathcal H\rightarrow\mathbb C^3$, $A_\zeta,\mathbb X(x,t),A:\mathcal H\rightarrow\mathcal H$, $B(x,t):\mathbb C^3\rightarrow\mathcal H$ 
and a $3\times 3$ matrix function $\gamma_*(x,t)$ satisfy
\[ \begin{array}{llllllll}
\dfrac{\partial}{\partial x} B &= & - (A B \sigma_2 + B \gamma)\sigma_1^{-1}& \dfrac{\partial}{\partial t} B & & = - (A \, B \widetilde\sigma_2 + B \widetilde\gamma)\widetilde\sigma_1^{-1}, \\
\dfrac{\partial}{\partial x} C &=& \sigma_1^{-1}( \gamma C - \sigma_2 C A_\zeta) & \dfrac{\partial}{\partial t} C  & & = \widetilde\sigma_1^{-1} (\widetilde\gamma C - \widetilde\sigma_2 C A_\zeta ), \\
\dfrac{\partial}{\partial x} \mathbb X & =&   B \sigma_2 C & \dfrac{\partial}{\partial t} \mathbb X &  & =  B \widetilde\sigma_2 C, \\
 0 &=& A \mathbb X + \mathbb X A_\zeta  + B \sigma_1 C, \\
\gamma_*& =&  \gamma + \sigma_2 C \mathbb X^{-1} B \sigma_1 - \sigma_1 C \mathbb X^{-1} B \sigma_2,
\end{array} \]
is a Boussinesq vessels. When the operators are bounded the function 
\[ q(x)=-\dfrac{3}{2} \dfrac{\partial^2}{\partial x^2} \ln (\tau(x,t)) = -\dfrac{3}{2} \dfrac{\partial^2}{\partial x^2} \ln \det(\mathbb X^{-1}(0,0)\mathbb X(x,t))
\]
satisfies the Boussinesq equation \begin{equation} \label{eq:Boussinesq}
q_{tt} = \dfrac{\partial^2}{\partial x^2} [3q_{xx} - 12 q^2 ]
\end{equation}
on the set $\Omega\subseteq\mathbb R^2$, where $\mathbb X(x,t)$ is invertible
(See \cite{bib:BoussVessels} for details). A similar to Main Theorem \ref{mthm} result can be proved for the Boussinesq equation \eqref{eq:Boussinesq} using ideas of this work.


\section{\label{sec:Back}Background on Krein space theory}
Let $(\mathcal H,\langle\cdot,\cdot\rangle)$ be a Hilbert space. Let $\mathbb X$ be a self-adjoint bounded operator on $\mathcal H$.
We define a sesquilinear form $[\cdot,\cdot]$ on $\mathcal H$ as $[u,v]=\langle \mathbb X u,v\rangle$. The map $[\cdot,\cdot]: \mathcal H\times\mathcal H\rightarrow\mathbb C^2$ possesses
the following properties
\begin{enumerate}
	\item $[\alpha x_1+\beta x_2,y] = [\alpha x_1,y] + [\beta x_2,y]$, $x_1,x_2,y\in\mathcal H$, $\alpha,\beta\in\mathbb C$,
	\item $[x,y]=\overline{[y,x]}$, $x,y\in\mathcal H$,
	\item $|[x,y]| \leq C \|x\|\|y\|$, $x,y\in\mathcal H$ with $C>0$.
\end{enumerate}
Here $\|\cdot\|$ denotes the Hilbert space norm of $\mathcal H$. Conversely, if there is a map $[\cdot,\cdot]$ with these three properties, there exists
a self-adjoint operator $\mathbb X:\mathcal H\rightarrow\mathcal H$ realizing it by $[x,y]=\langle\mathbb X x,y\rangle$.

Let $\mathcal K$ be equal to the Hilbert space $\mathcal H$ as a set, 
equipped with (indefinite) inner product: $(\mathcal K,[\cdot,\cdot])$. Then the pair
$(\mathcal K,[\cdot,\cdot])$ is called \textit{Krein space}. For any operator $T$
on $\mathcal K$ we denote by $T^*$\footnote{In the literature one usually denotes the adjoint of an operators $T$, with respect to $\langle\cdot,\cdot\rangle$ as $T^*$, and the adjoint with respect to
$[\cdot,\cdot]$ as $T^+$. Since we are dealing exclusively with the Krein-space adjoint, we will use $T^*$ for this notation and will rarely denote by
$T^{\langle*\rangle}$ the adjoint with respect to the Hilbert space $\mathcal H$.}
the unique operator satisfying $[Tu,v]=[u,T^*v]$ for all $u,v\in\mathcal K$. 
If we denote by $T^{<*>}$ the Hilbert space adjoint of $T$, and if $\mathbb X$ is invertible, then
\[ T^* = \mathbb X^{-1} T^{\langle*\rangle} \mathbb X.
\]
The space $\mathcal H$ admits the decomposition
\[ \mathcal H = \mathcal H^+ \oplus \mathcal H^-
\]
such that $[u,u]>0$ for all $x\in\mathcal H^+$ and $[u,u]<0$ for all $x\in\mathcal H^-$. Moreover, the spaces \linebreak
$(\mathcal H^+,[\cdot,\cdot]), (\mathcal H^-,-[\cdot,\cdot])$ are complete with respect to the norms 
$[\cdot,\cdot]$ and $-[\cdot,\cdot]$ respectively.

A typical example of a Krein space is as follows.
If $\rho(\mu)$ is a real function which is locally of bounded variation and $|\rho(\mu)|$ denotes its total variation, 
then the space $L^2(\rho)$, of all measurable functions $f$ such that $\int_{\mathbb R} |f(\mu)|^2 d|\rho(\mu)|<\infty$ and equipped with the indefinite inner product 
\[ [f,g] = \int_{\mathbb R} g^*(\mu) f(\mu) d\rho(\mu),
\]
is a Krein space. A typical example used in this work involves a matrix-valued measure. Let $d\bar\rho = \bbmatrix{d\rho_{11}&d\rho_{12}\\d\rho_{21}&d\rho_{22}}$ be a 
$2\times 2$ matrix of measures. A complex-valued, self-adjoint
\footnote{self-adjoint means $\rho_{11}=\rho^*_{11}$, $\rho_{12}=\rho^*_{21}$, $\rho_{22}=\rho^*_{22}$}, $2\times 2$ measure is called \textit{positive}, if
\[ \int_{\mathbb R} \bbmatrix{f_1^*(\mu) & f_2^*(\mu)} d\bar\rho(\mu) \bbmatrix{f_1(\mu)\\f_2(\mu)} \geq 0,
\]
provided the integral exists. Denote by $\mathcal R = supp(d\rho)$ - the support of the measure $d\rho$, then it is a matter of standard verifications that 
\[ \mathcal H = L^2(d\bar\rho) = \{ \bar f = \bbmatrix{f_1(\mu)\\f_2(\mu)}  \mid \|f\|^2 = \int_{\mathcal R} \bbmatrix{f_1^*(\mu) & f_2^*(\mu)} d\bar\rho(\mu) \bbmatrix{f_1(\mu)\\f_2(\mu)} < \infty \}
\]
is a Hilbert space. Suppose that $d\bar\rho=d\bar\rho_+-d\bar\rho_-$ for two positive measures $d\bar\rho_+$, $d\bar\rho_-$, creating
two Hilbert spaces of column-functions $\mathcal H_+$, $\mathcal H_-$ as above. We define
\[ \mathcal K = \{ \bar f(\mu) \mid \| \bar f \|_{\mathcal H_+} + \| \bar f \|_{\mathcal H_-} < \infty \},
\]
equipped with the indefinite inner product ($\mathcal R = supp(d\rho_+) \cup supp(d\rho_-)$)
\[ \begin{array}{lll}
[\bar f, \bar g] & =  \int_{\mathcal R} \bbmatrix{g_1^*(\mu) & g_2^*(\mu)} d\bar\rho(\mu) \bbmatrix{f_1(\mu)\\f_2(\mu)} \\
& = \int_{\mathcal R} \bbmatrix{g_1^*(\mu) & g_2^*(\mu)} d\bar\rho_+(\mu) \bbmatrix{f_1(\mu)\\f_2(\mu)}  -
\int_{\mathcal R} \bbmatrix{g_1^*(\mu) & g_2^*(\mu)} d\bar\rho_-(\mu) \bbmatrix{f_1(\mu)\\f_2(\mu)}.
\end{array} \]
The space of all bounded operators between Krein space is denoted by $L(\mathcal K_1,\mathcal K_2)$.
In this work we frequently use a 2 dimensional Hilbert space $\mathbb C^2$ for either $\mathcal K_1$ or $\mathcal K_2$. 
In this case we identify the sesquilinear form on $\mathbb C^2$
with the standard inner product of $\mathbb C^2$.

If $T:\mathcal K_1\rightarrow\mathcal K_2$ then its adjoint $T^*:\mathcal K_2\rightarrow\mathcal K_1$ is defined as the unique operator, satisfying
\[ [Tx,y ]_{\mathcal K_2} = [x,T^* y ]_{\mathcal K_1}.
\]
We present class of operators, which generate analytic semi groups. An operator $A:\mathcal K\rightarrow\mathcal K$
(usually instead of $\mathcal K$ Banach spaces are used) generates an analytic semigroup if there exists $w>0$ such that
$\Re\lambda > w$ is contained in the resolvent set of $A$ and there is $C>0$ such that
\[ \| (\lambda I - A)^{-1}\| \leq \dfrac{C}{|\lambda - w|}.
\]
The resolvent set of $A$ contains also the sector of the form
\[ Sec = \{ \lambda\in\mathbb C \mid |\arg(\lambda)-w) < \dfrac{\pi}{2} + \delta\}
\]
for some $\delta>0$. Such generators possess ``functional calculus'':
\begin{equation} \label{eq:fofgeneratorA}
f(A) = \dfrac{1}{2\pi i} \int\limits_\Gamma f(\lambda) (\lambda I - A)^{-1} d\lambda,
\end{equation}
where $f(\lambda)$ is analytic in $Sec$ and the curve $\Gamma$ goes from $e^{-i\Theta_0}\infty$
to $e^{i\Theta_0}\infty$ entirely inside of $Sec$ (with $\dfrac{\pi}{2} < \Theta_0 < \dfrac{\pi}{2} + \delta$).
For example, the analytic semigroup, generated by the operator $A$ is
\[ e^{Ax} = \dfrac{1}{2\pi i} \int\limits_\Gamma e^{\lambda x} (\lambda I - A)^{-1} d\lambda.
\]

We mention Hille-Yosida Theorem, characterizing generators of $C_0$ semigroups on $\mathbb R$, which is sufficient
for some of the theorems.
\begin{thm} \label{thm:HilleYosida} Let $A$ be a linear operator defined on a linear subspace $D(A)$ of the Banach space $\mathcal K$, $w$ be a real number, and $M > 0$. 
Then $A$ generates a strongly continuous semigroup, denoted as $e^{Ax}$, that satisfies $\| e^{Ax}\| \leq M e^{w x}$ if and only if
\begin{enumerate}
	\item $D(A)$ is dense in $\mathcal K$, and
	\item every real $\lambda > w$ belongs to the resolvent set of A and for such $\lambda$ and for all positive integers $n$:
	\[\|(\lambda I - A)^{-n}\| \leq \dfrac{M}{(\lambda-w)^n}.\]
\end{enumerate}
\end{thm}
A proof of this Theorem can be found in \cite[Theorem 3.4.1]{bib:Staffans}, \cite[Theorem II.3.5]{bib:EngelNagel}.

\section{Non-symmetric vessels}
Theory of operator nodes is presented in \cite{bib:Brodskii}. We use a generalization of this notion, involving unbounded operators.
This notion is used to study bounded operators $A = A_R + A_I = \dfrac{A+A^*}{2} + \dfrac{A-A^*}{2i}$, whose image part $A_I$ (or real part $A_R$) is small,
or more precisely is compact. We substitute this requirement by the existence of $A_\zeta:\mathcal H\rightarrow\mathcal H$,
such that $A+A_\zeta$ is 2-dimensional in a Krein space and these two operators have the same domain.

Many notions from the theory of nodes \cite{bib:Brodskii} can be applied to the notion of a node, presented in this work. We have not inserted these 
results primarily  for the lack of space, but also because of a different aim: we want to prove the existence of solutions for the KdV equation \eqref{eq:KdV}.

A prevessel is a node, for which some of the operators depend on $x\in\mathbb R$, and a vessel is an ``invertible'' prevessel (in the sense of Brodskii).
The use of unbounded operators requires a careful consideration of their domains. The axioms of a node, presented here, assume equations, which take this issue into account.

\subsection{Node, prevessel, vessel}
\begin{defn} A \textbf{node} is a collection of operators and spaces
\[ \mathfrak N = \bbmatrix{
C & A_\zeta, \mathbb X, A & B & \sigma_1\\
\mathbb C^2 & \mathcal K & \mathbb C^2}
\]
where $\mathcal K$ is a Krein space, $C:\mathcal K\rightarrow\mathbb C^2$,
$\mathbb X: \mathcal K\rightarrow\mathcal K$,  $B:\mathbb C^2\rightarrow\mathcal K$ are bounded operators,
$\sigma_1=\sigma_1^*$ - invertible $2\times 2$ matrix, 
$A, A_\zeta$ are generators of $C_0$ groups on $\mathcal K$ with identical dense domain $D(A)=D(A_\zeta)$. The operator $\mathbb X$
is assumed to satisfy $\mathbb X(D(A))\subseteq D(A)$. The operators of the node are subject to the \textbf{Lyapunov equation}
\begin{equation} \label{eq:Lyapunov}
A \mathbb X u + \mathbb X A_\zeta u + B \sigma_1 C u = 0, \quad \forall u\in D(A_\zeta)=D(A).
\end{equation}
If $\mathbb X$ is invertible, the \textbf{transfer function} of $\mathfrak N$ is
\begin{equation} \label{eq:S0realized}
S(\lambda) = I - C \mathbb X^{-1} (\lambda I - A)^{-1} B \sigma_1.
\end{equation}
The node $\mathfrak N$ is called \textbf{symmetric} if $A_\zeta=A^*$ and $C=B^*$.
\end{defn}
\Rems a function $S(\lambda)$, representable in the form \eqref{eq:S0realized} is called \textit{realized} \cite{bib:bgr}.
\textbf{2.} if $A_\zeta=A+T$ for a bounded operator $T$, then $D(A)=D(A_\zeta)$ ($A_\zeta$ is called a perturbation of $A$ in this case).
\textbf{3.} if $\mathbb X=I$, then the condition $\mathbb X(D(A))\subseteq D(A)$ holds. \textbf{4.} when the node is symmetric one can verify that
\[ S(\lambda) \sigma_1^{-1} S^*(-\bar\lambda) = \sigma_1^{-1}
\]
at all points of analyticity of $S$. \textbf{5.} for the unbounded operators $A,A_\zeta$ to be generators of $C_0$-groups, it necessary and sufficient
to demand that they satisfy the conditions of the Hille-Yosida Theorem \ref{thm:HilleYosida}.
Particularly, they must be closed, densely defined operators.

In the case $\mathbb X$ is invertible, we consider a stronger notion of a node as follows.
\begin{defn} A node $\mathfrak N$ is called \textbf{invertible}, if $\mathbb X$ is invertible and $\mathbb X^{-1}(D(A))\subseteq D(A)$.
\end{defn}
A simple chain of inclusions for an invertible node
\[ \mathbb X(D(A))\subseteq D(A) \Rightarrow D(A)\subseteq \mathbb X^{-1}(D(A)) \subseteq D(A),
\]
where the first inclusion comes from the node condition, and the last one from the invertible node condition, implies that
$\mathbb X^{-1}(D(A)) = D(A)$. Similarly, $\mathbb X(D(A)) = D(A)$. Moreover, taking $u=\mathbb X^{-1}u'$, where $u,u'\in D(A)$ and plugging it into the
Lyapunov equation \eqref{eq:Lyapunov}, we obtain that
\begin{equation} \label{eq:Lyapunov-1}
A_\zeta \mathbb X^{-1} u' + \mathbb X^{-1} A u' + \mathbb X^{-1}B \sigma_1 C \mathbb X^{-1}u' = 0, \quad \forall u'\in D(A),
\end{equation}
after multiplying by $\mathbb X^{-1}$ from the left. From the existence of this Lyapunov equation we obtain the following Lemma. 
\begin{lemma} If $\mathfrak N$ is an invertible node, then
\[ \mathfrak N_{-1} = \bbmatrix{
C\mathbb X^{-1} &  A, \mathbb X^{-1},A_\zeta & \mathbb X^{-1}B & \sigma_1\\
\mathbb C^2 & \mathcal K & \mathbb C^2}
\]
is also a node.
\end{lemma}
One could consider a similar notion of ``adjointable'' node, for which the adjoint of the Lyapunov equation \eqref{eq:Lyapunov} would define
a node, but we do not insert details here. Actually, there is a theory of construction of new such nodes from old ones, similarly to the theory
presented in \cite{bib:Brodskii, bib:bgr}. In the case $\mathbb X = I$ (the identity operator) we have a very well developed theory \cite{bib:BL} of 
(symmetric) nodes with $A_\zeta = A^*$, which has a finite dimensional real part: $A + A^* = -B\sigma_1B^*$. 

Finally, rewriting the Lyapunov equation \eqref{eq:Lyapunov-1}, of the invertible node as follows
\[ (-A_\zeta) \mathbb X^{-1} u' + \mathbb X^{-1} (-A) u' + \mathbb X^{-1}B \sigma_1^{-1} \sigma_1(-\sigma_1 C) \mathbb X^{-1}u' = 0, \quad \forall u'\in D(A)
\]
we arrive to the node
\[ \mathfrak N^{-1} = \bbmatrix{
-\sigma_1 C\mathbb X^{-1} &  -A, \mathbb X^{-1}, -A_\zeta &  \mathbb X^{-1}B\sigma_1^{-1} & \sigma_1\\
\mathbb C^2 & \mathcal K & \mathbb C^2}
\]
whose transfer unction
\[ S^{-1}(\lambda) = I + \sigma_1 C (\lambda I + A_\zeta)^{-1}\mathbb X^{-1}B
\]
is the inverse of the transfer function, defined in \eqref{eq:S0realized},
of the original invertible node $\mathfrak N$. This is a standard fact, related to Schur complements and can be found in \cite{bib:Brodskii, bib:bgr}.
\begin{defn} Class $\R(\sigma_1)$ consist of $2\times 2$ matrix-valued functions $S(\lambda)$ of the complex variable $\lambda$, which are transfer functions of invertible nodes.
The subclass $\U(\sigma_1)\subseteq\R(\sigma_1)$ consists of the transfer functions of symmetric, invertible nodes.
The Schur class $\SU(\sigma_1)\subseteq\U$ demand also that the inner space $\mathcal K$ is Hilbert and $\mathbb X>0$. The sub-classes of rational functions
in $\SU, \U, \R$ are denoted by $\RSU, \RU, \RR$ respectively.
\end{defn}
When $S(\lambda)$ is just analytic at infinity (hence $A$ must be bounded), 
there is a very well known theory of realizations developed in \cite{bib:bgr}. 
For analytic at infinity and symmetric, i.e. satisfying $S^*(-\bar\lambda) \sigma_1 S(\lambda) =  \sigma_1$,
functions there exists a good realization theory using Krein spaces,
developed in \cite{bib:KreinReal}\footnote{At the paper \cite{bib:KreinReal} a similar result is proved for functions symmetric with respect to the unit circle,
but it can be translated using Calley transform into $S^*(-\bar\lambda) \sigma_1 S(\lambda) =  \sigma_1$ and was done in \cite{bib:GenVessel, bib:SchurIEOT}}.
Such a realization is then translated into a function in $\U(\sigma_1)$.
The sub-classes $\U, \SU$ appear a lot in the literature and correspond to the symmetric case. We will not particularly
consider these two classes here and refer to \cite{bib:SchurIEOT}.

Equations, which arise in the theory of vessels involve differential equations
with unbounded operators. As a result, an operator satisfying such an equation must satisfy a relation
with the domain of the unbounded operator, which is presented in the next Definition.
\begin{defn} A bounded operator $B:\mathbb C^2\rightarrow\mathcal K$ is called $A$-\textbf{regular},
where $A: \mathcal K\rightarrow\mathcal K$ is linear, if $Be\in D(A)$ for all $e\in\mathbb C^2$.
\end{defn}
\begin{defn} \label{def:preVessel} Fix $2\times 2$ matrices $\sigma_2=\sigma_2^*$, $\gamma = -\gamma^*$.
The collection of operators and spaces 
\begin{equation} \label{eq:DefpreV}
\mathfrak{preV} = \bbmatrix{C(x) & A_\zeta, \mathbb X(x), A & B(x)&  \sigma_1,\sigma_2,\gamma \\
\mathbb C^2 & \mathcal{K} & \mathbb C^2}
\end{equation}
is called a (non-symmetric) \textbf{prevessel}, if the following conditions hold:
1. $\mathfrak{preV}$ is a node for all $x\in\mathbb R$, 2.
the operator $B(x)\sigma_2$ is $A$-regular, 3. $C(x), \mathbb X(x), B(x)$ are differentiable, bounded
operators, subject to the following conditions
\begin{eqnarray}
\label{eq:DB} \frac{\partial}{\partial x} B(x)  = - (A B \sigma_2 + B \gamma) \sigma_1^{-1}, \\
\label{eq:DC} \frac{\partial}{\partial x} C(x)u = \sigma_1^{-1}(- \sigma_2 C A_\zeta u + \gamma Cu) , \quad \forall u\in D(A_\zeta), \\
\label{eq:DX} \frac{\partial}{\partial x} \mathbb X =  B \sigma_2 C,
\end{eqnarray}
The prevessel $\mathfrak{preV}$ is called \textbf{symmetric} if
$A_\zeta=A^*$ and $C(x)=B^*(x)$ for all $x\in\mathbb R$.
\end{defn}
It turns out that the structure of a prevessel implies the Lyapunov equations \eqref{eq:Lyapunov}, \eqref{eq:Lyapunov-1} as the following
Lemma claims.
\begin{lemma}[\textbf{permanence of the Lyapunov equations}] \label{lemma:Redund}
Suppose that $B(x), C(x), \mathbb X(x)$ satisfy \eqref{eq:DB}, \eqref{eq:DC}, \eqref{eq:DX} respectively and $\mathbb X(x)(D(A_\zeta))\subseteq D(A)$
for all $x\in\mathbb R$. Then if the Lyapunov equation \eqref{eq:Lyapunov}
holds for a fixed $x_0$, then it holds for all $x$. In the case the operator $\mathbb X(x)$ is invertible and $B(x), C(x), \mathbb X(x)$ are part of an invertible node,
if \eqref{eq:Lyapunov-1} holds for a fixed $x_0$, then it holds for all $x$.
\end{lemma}
\noindent\textbf{Proof:} 
Let us differentiate the right hand side of the Lyapunov equation \eqref{eq:Lyapunov}:
\[ \begin{array}{lll}
\dfrac{d}{dx} [A \mathbb X u + \mathbb X A_\zeta u + B \sigma_1 C u ] = \\
= A B(x)\sigma_2C(x)  u + B(x)\sigma_2C(x)  A_\zeta u - A B(x) \sigma_2 C(x) u - B(x) \sigma_2 C(x) A_\zeta u \\
= 0.
\end{array} \]
The terms involving $\gamma$ are canceled, because $\gamma+\gamma^*=0$, by the assumption on it.
Thus it is a constant and the result follows. For the invertible node case, the condition \eqref{eq:Lyapunov-1} is a result of \eqref{eq:Lyapunov}. \qed

\begin{defn} \label{def:Vessel} The collection of operators, spaces and a set $\Omega\subseteq\mathbb R$ 
\begin{equation} \label{eq:DefVGen}
\mathfrak{V} = \bbmatrix{C(x) & A_\zeta, \mathbb X(x), A & B(x)&  \sigma_1,\sigma_2,\gamma,\gamma_*(x) \\
\mathbb C^2 & \mathcal{K} & \mathbb C^2 & \Omega}
\end{equation}
is called a (non-symmetric) \textbf{vessel}, if $\mathfrak{V}$ is a pre-vessel, $\mathbb X(x)$ is invertible on
$\Omega$, and $\mathfrak{V}$ is also an invertible node for all $x\in\Omega$.
The $2\times 2$ matrix-function $\gamma_*(x)$ satisfies the \textbf{linkage condition} on $\Omega$
\begin{equation} \label{eq:Linkage}
	\gamma_*  =  \gamma + \sigma_2 C \mathbb X^{-1} B \sigma_1 - \sigma_1 C \mathbb X^{-1} B \sigma_2.
\end{equation}
\end{defn}

The class of the transfer functions of vessels is defined as follows
\begin{defn} Class $\I=\I(\sigma_1,\sigma_2,\gamma;\Omega)$ consist of $2\times 2$ matrix-valued (transfer) functions $S(\lambda,x)$ of the complex variable $\lambda$ and $x\in\Omega\subseteq\mathbb R$, possessing the following representation:
\begin{equation} \label{eq:Srealized}
 S(\lambda,x) = I - C(x) \mathbb X^{-1}(x) (\lambda I - A)^{-1} B(x) \sigma_1,
\end{equation}
where the operators $C(x), \mathbb X(x), B(x)$ are part of a vessel $\mathfrak V$.
\end{defn}
Before we prove the B\" acklund transformation Theorem \ref{thm:Backlund} we present a technical lemma.
\begin{lemma} \label{lemma:X-1B} Let $\mathfrak{V}$ be a vessel. Then for all $u\in D(A)$
\begin{eqnarray}
\label{eq:DCX-1} \sigma_1 \dfrac{d}{dx} [C(x)\mathbb X^{-1}(x)] u = \sigma_2 C(x)\mathbb X^{-1}(x) A u + \gamma_*(x) C(x)\mathbb X^{-1}(x) u, \\
\label{eq:DX-1B} \dfrac{d}{dx} [\mathbb X^{-1}(x)B(x)] \sigma_1 = A_\zeta \mathbb X^{-1}(x) B(x)\sigma_2 - \mathbb X^{-1}(x) B(x) \gamma_*(x).
\end{eqnarray}
\end{lemma}
\Proof Consider \eqref{eq:DCX-1} first. We write under each equality the corresponding equation that is used to derive the next line:
\[ \begin{array}{lllllll}
\sigma_1 \dfrac{\partial}{\partial x} [C\mathbb X^{-1}] u & = \sigma_1 \dfrac{\partial}{\partial x} (C) \mathbb X^{-1} u +\sigma_1 C\dfrac{\partial}{\partial x} (\mathbb X^{-1}) u \\
 &  \quad \quad \text{\eqref{eq:DC}: } \quad \quad \sigma_1 \frac{\partial}{\partial x} C(x)u = - \sigma_2 C A_\zeta u + \gamma C u \\
 & \quad \quad \text{\eqref{eq:DX}:} \quad \quad \frac{\partial}{\partial x} \mathbb X =  B \sigma_2 C \\
 & = - \sigma_2 C A_\zeta \mathbb X^{-1} u + \gamma C \mathbb X^{-1} u  - \sigma_1 C\mathbb X^{-1} B \sigma_2 C \mathbb X^{-1} u \\
 & = -\sigma_2 C A_\zeta \mathbb X^{-1} u  + ( \gamma - \sigma_1 C\mathbb X^{-1} B \sigma_2) C \mathbb X^{-1} u \\
 & \quad \quad \text{\eqref{eq:Lyapunov-1}: }  \quad \quad  A_\zeta \mathbb X^{-1} u + \mathbb X^{-1} A u + \mathbb X^{-1}B \sigma_1 C \mathbb X^{-1}u = 0 \\
 & = \sigma_2 C \mathbb X^{-1} Au  + ( \gamma  + \sigma_2 C\mathbb X^{-1} B \sigma_1 - \sigma_1 C\mathbb X^{-1} B \sigma_2) C \mathbb X^{-1} u \\
 & \quad \quad \text{\eqref{eq:Linkage}: } \quad \quad  \gamma_*  =  \gamma + \sigma_2 C \mathbb X^{-1} B \sigma_1 - \sigma_1 C \mathbb X^{-1} B \sigma_2 \\
 & = \sigma_2 C\mathbb X^{-1} A u + \gamma_* C \mathbb X^{-1}u.
\end{array} \]
Notice that all equations of the vessel can be used, since we apply them to a vector $u$ from $D(A)$.

The equation \eqref{eq:DX-1B} is proved in exactly the same manner. \qed

Now we have all the ingredients of the following Theorem. This theorem has its origins at the work of M. Liv\c sic \cite{bib:Vortices} and was
proved for bounded operators in \cite{bib:GenVessel, bib:ENLS, bib:SchurIEOT}. Now we present a generalization of these results for the case of unbounded
operator $A$.
\begin{thm}[Vessel=B\" acklund transformation] \label{thm:Backlund}
Let $\mathfrak{V}$ be a vessel defined in \eqref{eq:DefVGen} and satisfying the conditions of Definition \ref{def:Vessel}. Fix
$\lambda\not\in\SPEC(A)$ and let $u(\lambda,x)$ be a solution of the input LDE 
\begin{equation} \label{eq:InCC}
 \lambda \sigma_2 u(\lambda, x) -
\sigma_1 \frac{\partial}{\partial x}u(\lambda,x) +
\gamma u(\lambda,x) = 0.
\end{equation} 
Then the function $y(\lambda,x)=S(\lambda,x)u(\lambda,x)$ satisfies the output LDE 
\begin{equation} \label{eq:OutCC}
\lambda \sigma_2 y(\lambda, x) - \sigma_1 \frac{\partial}{\partial x}y(\lambda,x) +
\gamma_*(x) y(\lambda,x) = 0.
\end{equation}
\end{thm}
\noindent\textbf{Proof:} Let us fix $\lambda\not\in\SPEC(A)$ and a solution $u(\lambda,x)$ of \eqref{eq:InCC}.
Then for $y(\lambda,x) = S(\lambda,x) u(\lambda,x)$ we calculate:
\[ \begin{array}{lll}
\sigma_1 \dfrac{d}{dx} y(\lambda,x)  & = \sigma_1 \dfrac{d}{dx} [(I - C(x)\mathbb X^{-1}(x) (\lambda I - A)^{-1} B(x)\sigma_1) u(\lambda,x)] =\\
& = \sigma_1 \dfrac{d}{dx} u(\lambda,x) - \sigma_1 \dfrac{d}{dx} [ C(x)\mathbb X^{-1}(x) (\lambda I - A)^{-1} B(x)\sigma_1 u(\lambda,x)] \\
& = (\sigma_2\lambda+\gamma)u(\lambda,x) - \sigma_1 \dfrac{d}{dx} [ C(x)\mathbb X^{-1}(x)] \, (\lambda I - A)^{-1} B(x)\sigma_1 u(\lambda,x) \\
& \quad \quad - \sigma_1C(x)\mathbb X^{-1}(x) (\lambda I - A)^{-1} \dfrac{d}{dx} [B(x)] \sigma_1 u(\lambda,x)  \\
&  \quad \quad \quad \quad - \sigma_1C(x)\mathbb X^{-1}(x) (\lambda I - A)^{-1} B(x) \sigma_1 \dfrac{d}{dx} u(\lambda,x) .
\end{array} \]
Using \eqref{eq:DCX-1}, \eqref{eq:DB} and \eqref{eq:InCC} it becomes (notice that $(\lambda I - A)^{-1} B(x)\sigma_1 u(\lambda,x)\in D(A)$)
\[ \begin{array}{llll}
\sigma_1 \dfrac{d}{dx} y(\lambda,x)   =
(\sigma_2\lambda+\gamma)u(\lambda,x) - \\
-  [\sigma_2 C(x)\mathbb X^{-1}(x) A + \gamma_*(x) C(x) \mathbb X^{-1}(x)] (\lambda I - A)^{-1} B(x)\sigma_1 u(\lambda,x) +\\
\quad \quad + \sigma_1C(x)\mathbb X^{-1}(x) (\lambda I - A)^{-1} [A B(x)\sigma_2 + B(x) \cancel{\gamma}] u(\lambda,x) \\
\quad \quad \quad \quad -\sigma_1C(x)\mathbb X^{-1}(x) (\lambda I - A)^{-1} B(x) (\sigma_2\lambda+\cancel\gamma) u(\lambda,x) =
\end{array} \]
Let us combine the last two terms and add $\pm \lambda I$ next to $A$:
\[ \begin{array}{llll}
= (\sigma_2\lambda + \gamma)u(\lambda,x) - \\
-  [\sigma_2 C(x)\mathbb X^{-1}(x) (A \pm \lambda I)+ \gamma_*(x) C(x) \mathbb X^{-1}(x)] (\lambda I - A)^{-1} B(x)\sigma_1 u(\lambda,x) +\\
\quad \quad + \sigma_1C(x)\mathbb X^{-1}(x) (\lambda I - A)^{-1} (A -\lambda I)B(x)\sigma_2 u(\lambda,x) = \\
= (\sigma_2\lambda + \gamma)u(\lambda,x) + \sigma_2 C(x)\mathbb X^{-1}(x)  B(x)\sigma_1 u(\lambda,x) - \\
- [\sigma_2 C(x)\mathbb X^{-1}(x) \lambda + \gamma_*(x) C(x) \mathbb X^{-1}(x)] (\lambda I - A)^{-1} B(x)\sigma_1 u(\lambda,x) - \\
\quad \quad - \sigma_1C(x)\mathbb X^{-1}(x) B(x)\sigma_2 u(\lambda,x) = \\
= (\sigma_2\lambda + \gamma + \sigma_2 C(x)\mathbb X^{-1}(x)  B(x)\sigma_1 -\sigma_1C(x)\mathbb X^{-1}(x) B(x)\sigma_2 )u(\lambda,x) - \\
- [\sigma_2 C(x)\mathbb X^{-1}(x) \lambda + \gamma_*(x) C(x) \mathbb X^{-1}(x)] (\lambda I - A)^{-1} B(x)\sigma_1 u(\lambda,x).
\end{array} \]
Using \eqref{eq:Linkage} and the definition of $S(\lambda,x)$ we obtain that
\[ \begin{array}{llll}
\sigma_1 \dfrac{d}{dx} y(\lambda,x) & = [\sigma_2 \lambda + \gamma_*(x)] u(\lambda,x) - \\
& \quad\quad [\sigma_2 \lambda  - \gamma_*(x)] C(x)\mathbb X^{-1}(x) (\lambda I - A)^{-1} B(x) \sigma_1 u(\lambda,x) = \\
&= [\sigma_2 \lambda + \gamma_*(x)] [I - C(x)\mathbb X^{-1}(x) (\lambda I - A)^{-1} B(x) \sigma_1] u(\lambda,x)=\\
&= (\sigma_2\lambda + \gamma_*(x)) S(\lambda,x) u(\lambda,x) = \\
&= (\sigma_2\lambda + \gamma_*(x)) y(\lambda,x).\qed
\end{array} \]
One of the corollaries \cite{bib:CoddLev} of this Theorem is that the function $S(\lambda,x)$ must satisfy the following differential equation
\begin{equation} \label{eq:DS}
\dfrac{\partial}{\partial x} S(\lambda,x) = \sigma_1^{-1} (\sigma_2\lambda+\gamma_*(x)) S(\lambda,x)- S(\lambda,x) \sigma_1^{-1} (\sigma_2\lambda+\gamma).
\end{equation}
Moreover, defining the fundamental solutions $\Phi(\lambda,x), \Phi_*(\lambda,x)$
\begin{eqnarray}
\label{eq:DefPhi} \lambda \sigma_2 \Phi(\lambda,x) -\sigma_1 \frac{\partial}{\partial x}\Phi(\lambda,x) + \gamma \Phi(\lambda,x) = 0, \quad \Phi(\lambda,0) = I, \\
\label{eq:DefPhi*} \lambda \sigma_2 \Phi_*(\lambda,x) - \sigma_1 \frac{\partial}{\partial x}\Phi_*(\lambda,x) +
\gamma_*(x) \Phi_*(\lambda,x) = 0, \quad \Phi_*(\lambda,0) = I,
\end{eqnarray}
we also obtain that
\begin{equation} \label{eq:SPhi*S0Phi}
S(\lambda,x) = \Phi_*(\lambda,x) S(\lambda,0) \Phi^{-1}(\lambda,x).
\end{equation}

\subsection{\label{sec:preVConcstr}Standard construction of a prevessel}
Now we present the standard construction of a prevessel $\mathfrak {preV}$ from a node $\mathfrak N_0$.
under assumption that the operators $A, A_\zeta$ are generators of analytic semi groups. In general, it is enough to demand that
$A, A_\zeta$ possess ``functional calculus''. Using formula \eqref{eq:fofgeneratorA} and the fundamental matrices $\Phi(\lambda,x), \Phi_*(\lambda,x)$,
defined in \eqref{eq:DefPhi}, \eqref{eq:DefPhi*} we make the following definition.
\begin{defn} Let 
\[ \mathfrak N_0 = \bbmatrix{
C_0 & A_\zeta, \mathbb X_0, A & B_0 & \sigma_1\\
\mathbb C^2 & \mathcal K & \mathbb C^2}
\]
be a node, such that $A, A_\zeta$ generate analytic semi groups (or possess ``functional calculus'') and $D(A)=D(A_\zeta)$. 
The \textbf{standard construction} of the operators $B(x), C(x), \mathbb X(x)$ from the node $\mathfrak N_0$ is as follows
\begin{eqnarray}
\label{eq:StConB} B(x) = \dfrac{1}{2\pi i} \int\limits_\Gamma (\lambda I - A)^{-1} B_0 \Phi^*(-\bar\lambda,x-x_0) d\lambda, \\
\label{eq:StConC} C(x) = \dfrac{1}{2\pi i} \int\limits_\Gamma \Phi(\lambda,x-x_0) C_0 (\lambda I + A_\zeta)^{-1} d\lambda, \\
\label{eq:StConX} \mathbb X(x) = \mathbb X_0 + \int\limits_{x_0}^x B(y)\sigma_2C(y) dy. 
\end{eqnarray}
\end{defn}
\begin{thm} \label{thm:STConpreVessel}The collection
\[ \mathfrak{preV} = \bbmatrix{C(x) & A_\zeta, \mathbb X(x), A & B(x)&  \sigma_1,\sigma_2,\gamma \\
\mathbb C^2 & \mathcal{K} & \mathbb C^2}
\]
defined by the standard construction from the node $\mathfrak N_0$ is a prevessel, coinciding with $\mathfrak N_0$ for $x=x_0$.
\end{thm}
\Proof The condition $B(x)\sigma_2$ is $A$-regular comes from the definition of $B(x)$. Indeed, for all
$\lambda,x$ $(\lambda I - A)^{-1} B_0 \Phi(\lambda,x-x_0) \in D(A)$. By the existence of the functional calculus, it follows that \eqref{eq:DB},
\eqref{eq:DC} hold. The equation \eqref{eq:DX} is immediate and the Lyapunov equation \eqref{eq:Lyapunov} follows from Lemma \ref{lemma:Redund}.
Finally, we have to show that $\mathbb X(x)(D(A)) \subseteq D(A)$. For each $u\in D(A)$
\[ \mathbb X(x) u = \mathbb X_0 u + \int_{x_0}^x B(y)\sigma_2 C(y) u dy.
\]
Here $\mathbb X_0 u\in D(A)$ by the assumptions on $\mathfrak N_0$. $B(y)\sigma_2\in D(A)$ by the $A$-regularity of $B(y)\sigma_2$. Moreover, since
for each $u\in D(A)$
\[ \dfrac{\partial}{\partial x} B(x)\sigma_1C(x) u = - A B(x)\sigma_2C(x)u - B(x)\sigma_2C(x)A_\zeta u \]
by integrating, we will obtain that
\[ \int_{x_0}^x A B(y)\sigma_2C(y)udy = B_0\sigma_1C_0 u - \dfrac{\partial}{\partial x} B(x)\sigma_1C(x) u + \int_{x_0}^x B(x)\sigma_2C(x)A_\zeta u dy
\]
exists. So, by the closeness of the operator $A$, we obtain that
\[ A \int_{x_0}^x B(y)\sigma_2 C(y) u dy = \int_{x_0}^x A B(y)\sigma_2C(y)udy
\]
exists and $\int_{x_0}^x B(y)\sigma_2 C(y) u dy \in D(A)$. 
\qed

If on an interval $\mathrm I$, including $x_0$ the operator $\mathbb X(x)$ is also invertible,
we can define $\gamma_*(x)$. In fact, the following Theorem holds.
\begin{thm}[local scattering] \label{thm:VFromS0} Suppose that $\mathfrak N_0$ is an invertible node, then there exists an interval $I$, including the given point $x_0$ and a vessel
$\mathfrak V$ on $\mathrm I$, such that at $x=x_0$ the vessel $\mathfrak V$ coincides with the node $\mathfrak N_0$.
\end{thm}
\Proof Since $\mathbb X_0$ is invertible, there exists a small interval $\mathrm I$, including $x_0$ on which the operator
\[ \mathbb X(x) = \mathbb X_0 + \int_{x_0}^x B(y)\sigma_2C(y)dy
\]
is invertible. As a result, we can define $\gamma_*(x)$ by the linkage condition \eqref{eq:Linkage}. In order to show that the collection
\eqref{eq:DefVGen}
\[
\mathfrak{V} = \bbmatrix{C(x) & A_\zeta, \mathbb X(x), A & B(x)&  \sigma_1,\sigma_2,\gamma,\gamma_*(x) \\
\mathbb C^2 & \mathcal{K} & \mathbb C^2 & \mathrm I }
\]
is a vessel, it is necessary and sufficient to show that $\mathfrak{V}$ is an invertible node for all $x\in\mathrm I$, for which in turn
we must show that $\mathbb X^{-1}(x)(D(A)) \subseteq D(A)$. Notice that
\[ \mathbb X^{-1}(x)u = \mathbb X_0^{-1}u - \int_{x_0}^x \mathbb X^{-1}(y)B(y)\sigma_2 C(y) \mathbb X^{-1}(y) u dy
\]
and
\begin{multline*}
 \dfrac{\partial}{\partial x} [\mathbb X^{-1}(x) B(x)\sigma_1C(x) \mathbb X^{-1}(x) u] = \\
A_\zeta \mathbb X^{-1}(x) B(x)\sigma_2C(x)\mathbb X^{-1}(x)u +
\mathbb X^{-1}(x)B(x)\sigma_2C(x)\mathbb X^{-1}(x)A,
\end{multline*}
following from \eqref{eq:DCX-1}, \eqref{eq:DX-1B}. As a result, we can use the same proof as for $\mathbb X(x)(D(A)) \subseteq D(A)$ in Theorem
\ref{thm:STConpreVessel}. \qed

The transfer function
\[ S(\lambda) = I - C_0\mathbb X_0^{-1}(\lambda I - A)^{-1} B_0\sigma_1
\]
can be considered as a ``scattering data'', because $\gamma_*(x)$ (a generalized
potential) is uniquely determined from $S(\lambda)$ by this construction. 
The uniqueness of $S(\lambda)$ for a given potential $\gamma_*(x)$ is false. For example, notice that multiplying the given intial value
$S(\lambda)$ by arbitrary scalar function $a(\lambda)$, bounded at infinity, with limit 1 there, we will obtain that the two functions
\[ \Phi_*(\lambda,x) S(\lambda)\Phi^{-1}(\lambda,x),\quad a(\lambda)\Phi_*(\lambda,x) S(\lambda)\Phi^{-1}(\lambda,x)
\]
correspond to the same $\gamma_*(x)$. They can be obtained by applying the standard construction to $S(\lambda)$ and to $a(\lambda)S(\lambda)$.

A weaker form of the uniqueness is presented in the next Lemma.
We emphasize that a similar Lemma was proved in 
the Sturm-Liouville case in \cite{bib:GenVessel} and analogous result exists in \cite{bib:FaddeyevII} 
for purely continuous spectrum.
\begin{lemma}\label{lemma:uniquegamma*} Suppose that two functions $S(\lambda,x)$, $\widetilde S(\lambda,x)$ are in class $\I(\sigma_1, \sigma_2,\gamma;\Omega)$,
possessing the same initial value
\[ S(\lambda,0) = \widetilde S(\lambda,0)
\]
and are bounded at a neighborhood of infinity, with a limit value $I$ there.
Then the corresponding outer potentials are equal on $\Omega$:
\[ \gamma_*(x) = \widetilde \gamma_*(x).
\]
\end{lemma}
\Proof Suppose that
\[ S(\lambda,x) = \Phi_*(\lambda,x) S(\lambda,0) \Phi^{-1}(\lambda,x), \quad \widetilde S(\lambda,x) = \widetilde \Phi_*(\lambda,x) S(\lambda,0) \Phi^{-1}(\lambda,x),
\]
as in \eqref{eq:SPhi*S0Phi}. Then
\[ \widetilde S^{-1}(\lambda,x) S(\lambda,x) = \widetilde \Phi_*(\lambda,x)\Phi^{-1}_*(\lambda,x)
\]
is entire (the singularities appear in $S(\lambda,0) = \widetilde S(\lambda,0)$ only and are canceled) and equal to $I$ (- the identity matrix) at infinity.
By a Liouville theorem, it is a constant function, namely $I$. So $\widetilde \Phi_*(\lambda,x)\Phi^{-1}_*(\lambda,x) = I$ or
\[ \widetilde \Phi_*(\lambda,x) = \Phi_*(\lambda,x).
\]
If we differentiate this, we obtain that $\widetilde \gamma_*(x)=\gamma_*(x)$ on $\Omega$.
\qed

\subsection{Tau function of a prevessel}
Existence of the vessel and its transfer function relies on the invertability of the operator $\mathbb X(x)$. In order to investigate the existence of the inverse
for $\mathbb X(x)$ notice that from \eqref{eq:DX}
\[ \mathbb X(x) = \mathbb X_0 + \int\limits_0^x B(y) \sigma_2 C(y) dy
\]
it follows that
\[ \mathbb X_0^{-1} \mathbb X(x) = I + \mathbb X^{-1}_0 \int\limits_{x_0}^x B(y) \sigma_2 C(y) dy.
\]
Since $\sigma_2$ has finite rank, 
this expression is of the form $I + T$, for a trace class operator $T$ and since 
$\mathbb X_0$ is an invertible operator, there exists a non trivial interval (of length at least $\dfrac{1}{\|\mathbb X_0^{-1}\|}$) on which $\mathbb X(x)$ and $\tau(x)$ are defined. Recall \cite{bib:GKintro} that a function $F(x)$ from (a, b) into the group G (the set of bounded invertible operators on H of the form I + T, for
a trace-class operator $T$) is said to be differentiable if $F(x) -I$ is \textit{differentiable} as a map into the trace-class operators. In our case,
\[ \dfrac{d}{dx} (\mathbb X_0^{-1}\mathbb X(x)) = 
\mathbb X_0^{-1} \dfrac{d}{dx} \mathbb X(x) =
\mathbb X_0^{-1} B(x)\sigma_2C(x)
\]
exists in trace-class norm. This leads us to the following
\begin{defn} \label{def:Tau} For a given prevessel $\mathfrak{preV}$ \eqref{eq:DefpreV} the tau function $\tau(x)$ is defined as
\begin{equation} \label{eq:Tau} \tau = \det (\mathbb X_0^{-1} \mathbb X(x)).
\end{equation}
\end{defn}
Israel Gohberg and Mark Krein \cite[formula 1.14 on p. 163]{bib:GKintro}
proved that if $\mathbb X_0^{-1}\mathbb X(x)$ is a differentiable function
into G, then $\tau(x) = \SP (\mathbb X_0^{-1}\mathbb X(x))$
\footnote{$\SP$ - stands for the trace in the infinite dimensional space.} is a differentiable map into $\mathbb C^*$ with
\begin{multline} \label{eq:GKform}
\dfrac{\tau'}{\tau}  = \SP (\big(\mathbb X_0^{-1} \mathbb X(x)\big)^{-1} 
\dfrac{d}{dx} \big(\mathbb X_0^{-1} \mathbb X(x)\big)) = \SP (\mathbb X(x)' \mathbb X^{-1}(x)) = \\
= \SP (B(x)\sigma_2 C(x) \mathbb X^{-1}(x)) =
\TR (\sigma_2 C(x) \mathbb X^{-1}(x)B(x)).
\end{multline}
Differentiating this expression, we obtain that
\[ (\dfrac{\tau'}{\tau})' = \dfrac{\tau''}{\tau} - (\dfrac{\tau'}{\tau})^2 =
\dfrac{d}{dx} \TR (\sigma_2 C(x) \mathbb X^{-1}(x)B(x)).
\]
Using vessel conditions, since $B(x)$, $\mathbb X^{-1}(x)$ are differentiable bounded operators in the case
$A\mathbb X^{-1}(x)B(x)$ exists, or in the case it is canceled (SL case) we obtain that
\begin{equation} \label{eq:tau'tau'}
 \begin{array}{llll}
\dfrac{d}{dx} \TR (\sigma_2 C(x) \mathbb X^{-1}(x)B(x)) = \\
=  \TR (\sigma_2 \sigma_1^{-1}(-\sigma_2 C(x)A_\zeta - \gamma^* C(x))\mathbb X^{-1}(x)B(x)) - \\
\quad - \TR(\sigma_2 C(x) \mathbb X^{-1}(x)B(x) \sigma_2 C(x) \mathbb X^{-1}(x) B(x))
+ \TR(\sigma_2 C(x) \mathbb X^{-1}(x) (-A B(x)\sigma_2 - B(x)\gamma) \sigma_1^{-1}) = \\
= \TR (\sigma_2 \sigma_1^{-1} \sigma_2 C(x)[ -A_\zeta \mathbb X^{-1} - \mathbb X^{-1} A]B(x)) - 
\TR( [\sigma_2 \sigma_1^{-1}\gamma^* + \gamma\sigma_1^{-1}\sigma_2] C(x)\mathbb X^{-1}(x)B(x)) - \\
\quad - \TR(\sigma_2 C(x) \mathbb X^{-1}(x)B(x) \sigma_2 C(x) \mathbb X^{-1}(x) B(x))= \\
= \TR (\sigma_2 \sigma_1^{-1} \sigma_2 C(x) \mathbb X^{-1}B(x)\sigma_1 C(x) \mathbb X^{-1} B(x)) - 
\TR( [\sigma_2 \sigma_1^{-1}\gamma^* + \gamma\sigma_1^{-1}\sigma_2] C(x)\mathbb X^{-1}(x)B(x)) - \\
\quad - \TR(\sigma_2 C(x) \mathbb X^{-1}(x)B(x) \sigma_2 C(x) \mathbb X^{-1}(x) B(x)).
\end{array} \end{equation}
\begin{thm} Suppose that $\mathfrak{preV}$ \eqref{eq:DefpreV} is a prevessel. Define an open set
\[ \Omega = \{ x \mid \tau(x) \neq 0 \}.
\]
Then the prevessel $\mathfrak{preV}$ is a vessel on $\Omega$.
\end{thm}
\Proof for each $x_0$ in which $\tau(x_0)\neq 0$, the operator $\mathbb X(x_0)$ is invertible. Then
there exists a closed interval $\mathrm I_{x_0}$, including $x_0$ on which  $\mathfrak{preV}$ defines a vessel
by Theorem \ref{thm:VFromS0}. Then $\Omega = \cup_{x_0} \mathrm I_{x_0}$ and it finishes the proof. \qed

\subsection{Moments and their properties}
If the function $S(\lambda)$ is analytic at the neighborhood of infinity, one can consider its Taylor series
\begin{equation} \label{eq:DefHn}
S(\lambda) = 
I - C\mathbb X^{-1}(\lambda I - A)^{-1}B\sigma_1 = I - \sum_{n=0}^\infty \dfrac{ H_n}{\lambda^{n+1}}\sigma_1.
\end{equation}
But in the general case this expansion may not converge. Still the following Definition can be in force.
\begin{defn}\label{defn:Hn} The \textbf{$n$-th moment} $H_n$ of a vessel $\mathfrak{V}$ is:
\[ H_n = C \mathbb X^{-1} A^n B,
\]
provided that the image of $B$ is in $D(A^n)$. The moment is defined to be infinity, otherwise.
\end{defn}
Moments will play a crucial role in the future research and we will the following defining property for them.
\begin{thm} Suppose that the moments $H_0,\ldots, H_{n+1}$ are finite and differentiable, then
\begin{equation}
\label{eq:DHn} 	(H_n)'_x = \sigma_1^{-1}\sigma_2H_{n+1} - H_{n+1} \sigma_2\sigma_1^{-1} + \sigma_1^{-1} \gamma_* H_n - H_n \gamma\sigma_1^{-1}.
\end{equation}
\end{thm}
\Proof In the regular case, when all the operators are bounded, \eqref{eq:DHn} follows from the differential equation \eqref{eq:DS}. 
In the general case, it follows from \eqref{eq:DB}, \eqref{eq:DX} and \eqref{eq:Linkage}.
\qed

\begin{thm}[Uniqueness of the moments] \label{thm:momentUnique} Suppose that two sequences of moments $H_n(x)$ and $\widetilde H_n(x)$ are finite, differentiable 
and satisfy \eqref{eq:Hi+1HiFs} with analytic $\gamma_*(x)$ and $\widetilde \gamma_*(x)$ respectively. Then from
\[ H_n(0) = \widetilde H_n(0), \quad \forall n=0,1,2,\ldots
\]
it follows that $\gamma_*(x)=\widetilde \gamma_*(x)$. If the infinite system of equations \eqref{eq:Hi+1HiFs} has
a unique sequence of solutions $H_n(x)$ for a given $\gamma_*(x)$ and initial values $H_n(0)$ then $H_n(x)=\widetilde H_n(x)$.
\end{thm}
\Proof Let us show by the induction that $H_0^{(n)}(0) = (\widetilde H)_0^{(n)}(0)$ for all $n=0,1,2,\ldots$.
And since these two moments are analytic, the result will follow from the uniqueness of the Taylor series.
For $n=0$, $H_0(0)=\widetilde H_0(0)$ and the basis of the induction follows. Then from \eqref{eq:DHn} it follows that
\[
H_0^{(1)}(x) = \sigma_1^{-1}\sigma_2 H_1 - H_1 \sigma_2\sigma_1^{-1} + \sigma_1^{-1} \gamma_* H_0 - H_0 \gamma\sigma_1^{-1}.
\]
Differentiating again, using \eqref{eq:DHn} for $n=0,1$ and the Linkage condition \eqref{eq:Linkage}, we will obtain that
\begin{multline*}
 H_0^{(1)}(x) = \sigma_1^{-1}\sigma_2 H_1 - H_1 \sigma_2\sigma_1^{-1} + \sigma_1^{-1} \gamma_* H_0 - H_0 \gamma\sigma_1^{-1} =\\
 = \sigma_1^{-1}\sigma_2 H_1 - H_1 \sigma_2\sigma_1^{-1} + \sigma_1^{-1} (\gamma + \sigma_2H_0\sigma_1 - \sigma_1H_0\sigma_2)H_0 - H_0 \gamma\sigma_1^{-1}
 = P_2(H_0(x),H_1(x),H_2(x)) 
\end{multline*}
for a non-commutative polynomial $P_2$ with constant coefficients (depending on $\sigma_1,\sigma_1,\gamma$). This shows that a simple
induction results in
\[ H_0^{(n)}(x) = P_n(H_0(x),H_1(x),\ldots,H_n(x),H_{n+1}(x))
\]
for a non-commutative polynomial $P_n$ with constant matrix-coefficients.
As a result, plugging here $x=0$ and using the condition $H_n(0) = \widetilde H_n(0)$
\begin{multline*}
 H_0^{(n)}(0) = P_n(H_0(0),H_1(0),\ldots,H_n(0),H_{n+1}(0)) = \\
= P_n(\widetilde H_0(0),\widetilde H_1(0),\ldots,\widetilde H_n(0),\widetilde H_{n+1}(0))=
(\widetilde H)_0^{(n)}(0).
\end{multline*}
From here it follows that $H_0(x)=\widetilde H_0(x)$ and hence by the linkage condition \eqref{eq:Linkage} $\gamma_*(x)=\widetilde \gamma_*(x)$. Then
the last statement of the Theorem follows from the uniqueness of solutions. \qed

\section{Sturm-Liouville vessels}
In the special case of SL vessel parameters, we obtain that equations \eqref{eq:InCC}, \eqref{eq:OutCC} are equivalent to \eqref{eq:SL}. Let us explain it
in more details.
\begin{defn} \label{def:SLparam}
The Sturm Liouville (SL) vessel parameters are defined as follows
\[ \sigma_1 = \bbmatrix{0 & 1 \\ 1 & 0}, \quad
\sigma_2 = \bbmatrix{1 & 0 \\ 0 & 0}, \quad
\gamma = \bbmatrix{0 & 0 \\ 0 & i}.
\]
\end{defn}
Suppose that we are given a SL vessel $\mathfrak{V}$, in other words, $\mathfrak{V}$ \eqref{eq:DefVGen} is defined
for the SL vessel parameters. Denote a differentiable $2\times 2$ matrix function 
$H_0 = B^*(x) \mathbb X^{-1}(x) B(x) = \bbmatrix{a&b\\c&d}$. 
Then from \eqref{eq:GKform} it follows that
$\dfrac{\tau'}{\tau}  = \TR (\sigma_2 B^*(x) \mathbb X^{-1}(x)B(x)) = a$ is the $1,1$ entry of $H_0$.
Using simple calculations it follows that
\[\sigma_2\sigma_1^{-1}\sigma_2=0, \quad \TR (\sigma_2 H_0\sigma_2H_0) = (\TR(\sigma_2H_0))^2, \quad
\sigma_2 \sigma_1^{-1}\gamma^* + \gamma\sigma_1^{-1}\sigma_2 = \bbmatrix{0&-i\\i&0} , \] 
and we obtain from \eqref{eq:tau'tau'} that
\[ \dfrac{\tau''}{\tau} = - \TR(\bbmatrix{0&-i\\i&0} H_0) = i (c-b).
\]
Notice that the terms involving operator $A$ are canceled. Moreover, we obtain that 
\[ \gamma_*(x) = \gamma + \sigma_2 H_0 \sigma_1 - \sigma_1 H_0 \sigma_2 =
\bbmatrix{b-c & a \\ -a & i} = \bbmatrix{i \dfrac{\tau''}{\tau} & \dfrac{\tau'}{\tau} \\ -\dfrac{\tau'}{\tau} & i}.
\]
Thus we obtain the following lemma (appearing already in \cite[Proposition 3.2 ]{bib:SLVessels})
\begin{lemma} \label{lemma:gamma*}For SL vessel parameters, the following formula for $\gamma_*(x)$ holds
\[ \gamma_*(x) = \bbmatrix{i \dfrac{\tau''}{\tau} & \dfrac{\tau'}{\tau} \\ -\dfrac{\tau'}{\tau} & i}.
\]
\end{lemma}

Analogously to \cite[Section 3.1.1]{bib:SLVessels}, simple calculations show that denoting $u(\lambda,x) = \bbmatrix{u_1(\lambda,x)\\ u_2(\lambda,x)}$ we shall obtain that the input compatibility
condition \eqref{eq:InCC} is equivalent to
\[ \left\{ \begin{array}{lll}
-\frac{\partial^2}{\partial x^2} u_1(\lambda,x) = -i\lambda u_1(\lambda,x), \\
u_2(\lambda,x) = - i \frac{\partial}{\partial x} u_1(\lambda,x).
\end{array}\right.
\]
The output $y(\lambda,x) = \bbmatrix{y_1(\lambda,x)\\y_2(\lambda,x)} = S(\lambda,x) u(\lambda,x)$ satisfies the output
equation \eqref{eq:OutCC}, which is equivalent to
\begin{equation} \label{eq:ysimpleODEs}
 \left\{ \begin{array}{lll}
-\frac{\partial^2}{\partial x^2} y_1(\lambda,x) - 2 \dfrac{d^2}{dx^2} [\ln(\tau(x))]  y_1(\lambda,x) = -i\lambda y_1(\lambda,x), \\
y_2(\lambda,x) = - i [ \frac{\partial}{\partial x} y_1(\lambda,x) + \dfrac{\tau'}{\tau}y_1(\lambda,x)].
\end{array}\right.
\end{equation}
Observing the first coordinates $u_1(\lambda,x), y_1(\lambda,x)$ of the vector-functions $u(\lambda,x),y(\lambda,x)$ we can see that multiplication by
$S(\lambda,x)$ maps solution of the trivial SL equation (i.e. $q(x)=0$) to solutions of the more complicated one,
defined by the potential 
\begin{equation} \label{eq:qbeta}
q(x) = - 2  \dfrac{d^2}{dx^2} [\ln(\tau(x))],
\end{equation}
which can be considered as an analogue of the scattering theory.
\subsection{\label{sec:S0Constr}Construction of a realized function, possessing given Moments}
Notice that the formula for the first moment is
\[ H_0 = \bbmatrix{ \dfrac{\tau'}{\tau} & b \\ c & h_0^{22}}, \quad i(b-c) = \dfrac{\tau''}{\tau}.
\]
In the case $c=-b$, we obtain the following form, which will be used further in the text:
\begin{equation}\label{eq:H0Mat}
H_0 = \bbmatrix{ \dfrac{\tau'}{\tau} & -\dfrac{i\tau''}{2\tau} \\ \dfrac{i\tau''}{2\tau} & h_0^{22}}.
\end{equation}

Studying equation \eqref{eq:DHn} one can obtain a formula for the moment $H_{n+1} = \bbmatrix{H_{n+1}^{11}&H_{n+1}^{12}\\H_{n+1}^{21}&H_{n+1}^{22}}$ in terms 
of $H_n = \bbmatrix{H_n^{11}&H_n^{12}\\H_n^{21}&H_n^{22}}$ as follows (here $\beta = -\dfrac{\tau'}{\tau}, \pi_{11} = \beta'-\beta^2$):
\begin{equation} \label{eq:Hi+1HiFs}
\left\{ \begin{array}{llll}
H_{n+1}^{11} & = iH_n^{22} -\dfrac{d}{dx} H_n^{12} + \beta H_n^{12} , \\
H_n^{12}-H_n^{21} & = i( \frac{d}{dx} H_{n+1}^{11} - \beta H_{n+1}^{11}) , \\
\frac{d}{dx} (H_n^{12} + H_n^{21}) & = -  i \pi_{11} H_n^{11} + \beta (H_n^{12} - H_n^{21}), \\
2i \frac{d}{dx} H_n^{22} & = \frac{d^2}{dx^2} H_n^{12} - 2 \beta \frac{d}{dx}H_n^{12}.
\end{array} \right. \end{equation}
Indeed, plugging the SL vessel parameters (Definition \ref{def:SLparam}) into \eqref{eq:DHn} we find that
\[\begin{array}{llll}
\frac{d}{dx} \bbmatrix{H_n^{11}&H_n^{12}\\H_n^{21}&H_n^{22}} & = \sigma_1^{-1}\sigma_2H_{n+1} - H_{n+1} \sigma_2\sigma_1^{-1} + \sigma_1^{-1} \gamma_* H_n - H_n \gamma\sigma_1^{-1}\\
& = \bbmatrix{0&-H_{n+1}^{11}\\H_{n+1}^{11}&H_{n+1}^{12}-H_{n+1}^{21}} + 
\bbmatrix{\beta H_n^{11} + i\beta(H_n^{21}-H_n^{12}) & \beta H_n^{12} + i H_n^{22}\\
-i\pi_{11}H_n^{11}-\beta H_n^{21}-H_n^{22} & -i\pi_{11}H_n^{12}-\beta H_n^{22}}
\end{array} \]
Expressing $H_{n+1}^{12}$ from $1,2$ and $2,1$ entries of this equality, we obtain that
\[ H_{n+1}^{11} = - \frac{d}{dx} H_n^{12} + \beta H_n^{12} + i H_n^{22} = \frac{d}{dx} H_n^{21} - (-i\pi_{11}H_n^{11}-\beta H_n^{21}-i H_n^{22})
\]
The first equality here is identical to the first equality in \eqref{eq:Hi+1HiFs}, and the second equality is identical to the third line of \eqref{eq:Hi+1HiFs}. 
In a similar manner one can derive the other equations (see \cite{bib:SchurVessels} for more details).

From the formulas \eqref{eq:Hi+1HiFs} it follows that we can construct the moments $H_n$ with a special pattern. Namely, one can choose the initial conditions for 
$H_n^{12} + H_n^{21}$ and $H_n^{22}$ (which are not assigned in \eqref{eq:Hi+1HiFs}) so that the following Lemma holds.
\begin{lemma} There exists a choice of initial conditions such that
\begin{equation} \label{eq:HnSpcial}
 H_n = i^n \bbmatrix{r_n^{11}& i b_n^{12} \\ -i b_n^{12}& d_n^{22}}
\end{equation}
with real-valued function $r_n^{11},b_n^{12},d_n^{22}$. More precisely, the conditions at $0$ must be chosen so that
\[ i^n H_n^{22}(0)\in\mathbb R, \quad H_n^{12} + H_n^{21} = 0.
\]
\end{lemma}
\Proof  Using induction, it is necessary to choose $r=0$, $h_0^{22}\in\mathbb R$ in formula \eqref{eq:H0Mat} and the first moment
$H_0$ will satisfy the required condition. Suppose that $H_n$ is of the form stated at the theorem. Then from \eqref{eq:Hi+1HiFs}
it follows that $H_{n+1}^{11} = i^{n+1} (d_n^{22}-\dfrac{d}{dx} r_n^{11} + \beta r_n^{11}) = i^{n+1} r_{n+1}^{11}$ and is of the required form.
Similarly, the other formulas of \eqref{eq:Hi+1HiFs} produce the required result. Notice that the initial conditions for 
$H_{n+1}$ must be chosen so that $i^{n+1}H_{n+1}^{22}(0) \in\mathbb R$ and $H_{n+1}^{12}(0) + H_{n+1}^{21}(0) = 0$. \qed

Suppose that we are given moments $H_n$, realized in the form \eqref{eq:HnSpcial}.
We would like to construct a function $S(\lambda)$ with these moments, analytic at $\mathbb C\backslash i\mathbb R$ and possessing the following realization form
\begin{equation} \label{eq:SindrhoForm}
S(\lambda) = I - \int\limits_{\mathbb R} \dfrac{1}{\lambda-i\mu} d\bar\rho(\mu)\sigma_1,
\end{equation}
where $d\bar\rho=d\bar\rho_+-d\bar\rho_-$ is a $2\times 2$ matrix measure, which is a difference of two positive matrix-measures $d\bar\rho_+$, 
$d\bar\rho_-$. Moreover, we want the measures $\bar\rho$, $d\bar\rho_+$, 
$d\bar\rho_-$ to be analytic, namely, to satisfy the following
\begin{defn}
A measure $d\bar\rho$ is called \textbf{analytic}, if $\int\limits_{\mathbb R} \mu^n d\bar\rho(\mu)$ is finite for each $n=0,1,2,\ldots$.
For each function $S(\lambda)$, realized in the form \eqref{eq:SindrhoForm} with an analytic measure $\bar\rho$, we define its $n$-th moment
as $H^S_n = \int\limits_{\mathbb R} (i\mu)^n d\bar\rho$.
\end{defn}
\begin{thm} \label{thm:S0FromHn}
Given a sequence of moments $H_n$ of the form \eqref{eq:HnSpcial}, there exists a function $S(\lambda)$,
realized in form \eqref{eq:SindrhoForm} with an analytic measure $d\bar\rho=\bbmatrix{d\rho_{11}&id\rho_{12}\\-id\rho_{12} & d\rho_{22}}$ on $\mathbb R$, 
with real signed measures $\rho_{11}, \rho_{12}, \rho_{22}$ whose moments
coincide with the given ones. Namely, it holds that $H^S_n = H_n$.
\end{thm}
\Rem this Theorem is best illustrated if $d\bar\rho$ has a bounded support. Taking $\lambda$, satisfying $|\lambda| > \sup|\operatorname{supp}(d\rho)|$, there exists Taylor
expansion $\dfrac{1}{\lambda-i\mu} = \sum\limits_{n=0}^\infty \dfrac{(i\mu)^n}{\lambda^{n+1}}$ and the moments $H_n^S$ are just the Taylor coefficients. This also gives an idea of the
converse construction, when given moments $H_n$ have exponential growth: $\| H_n \| \leq k C^n$.
Define the function as $S(\lambda) = I - \sum\limits_{n=0}^\infty \dfrac{H_n}{\lambda^{n+1}}\sigma_1$. The fact that this function can be extended to $\mathbb C\backslash iR$ with a 
realization in the form \eqref{eq:SindrhoForm} follows from this theorem.

\Proof We use the Hamburger problem, which constructs a positive Borel measure $\rho$ for a given set of real numbers $m_0,m_1,m_2,\ldots$.
More precisely, there exists a Borel measure $\rho$ satisfying $\int\limits_{\mathbb R} \mu^n d\rho(\mu) = m_n$ if and only if the Hankel matrices
\[ M_n = \bbmatrix{
m_0&m_1&\ldots&m_n\\
m_1&m_2&\ldots&m_{n+1}\\
\vdots&\vdots&\ddots&\vdots\\
m_n&m_{n+1}&\ldots& m_{2n}}
\]
are positive. 

In order to prove our theorem, notice that it is enough to prove the scalar case,
because the $2\times 2$ case consists of four problems for each entry:
\[ H_n^S = i^n\int\limits_{\mathbb R} \mu^n d\bar\rho(\mu) = i^n \int\limits_{\mathbb R}\mu^n \bbmatrix{d\rho_{11}&id\rho_{12}\\-id\rho_{12}&d\rho_{22}} = 
i^n \bbmatrix{r_n^{11}& i b_n^{12} \\ -i b_n^{12}& d_n^{22}} = H_n,
\]
or equivalently
\[ \int\limits_{\mathbb R} \mu^n d\rho_{11} = r_n^{11}, \quad 
\int\limits_{\mathbb R} \mu^n d\rho_{12} = b_n^{12}, \quad 
\int\limits_{\mathbb R} \mu^n d\rho_{22} = d_n^{22}.
\]
The scalar problem is solvable as follows. Suppose that we are given a set of real numbers $m_0,m_1,m_2,\ldots$. Let $M_n$ be the corresponding Hankel matrices. We are going to construct two additional sequences $v_n$, $u_n$ which
satisfy the conditions of the Hamburger Theorem and are such that $m_n = v_n - u_n$. For $m_0$, it is
immediate that $m_0=v_0-u_0$, for some positive $v_0\neq 0$ and $u_0$. 
Suppose by induction, that we have constructed $v_1, \ldots, v_{2n}$, 
$u_1, \ldots, u_{2n}$ and matrices $V_n$, $U_n$, similar to the construction of $M_n$. Suppose also by the induction 
that $\det V_n > 0$. Let $v_{2n+1}, v_{2(n+1)}$ be parameters and construct next
\[
M_{n+1} = \bbmatrix{M_n & \bar m\\\bar m^t & m_{2(n+1)}}, \quad
V_{n+1} = \bbmatrix{V_n & \bar v\\\bar v^t & v_{2(n+1)}}, \quad
U_{n+1} = V_{n+1}-M_{n+1},
\]
where $\bar m =\bbmatrix{m_{n+1}\\m_{n+2}\\\vdots\\m_{2n+1}}$, $\bar v =\bbmatrix{v_{n+1}\\v_{n+2}\\\vdots\\v_{2n+1}}$.
Using the principal minors criteria for the positivity of a matrix, all the principal minors of $V_{n+1}$ are those of
$V_n$, except for the last one:
\[
\det(V_{n+1}) = v_{2(n+1)} \det V_n + C_1.
\]
The last formula is obtained by expansion along the last column of $V_{n+1}$. Taking 
$v_{2(n+1)} > -\dfrac{C_1}{\det(V_n)}$, we obtain that $V_{n+1}$ is positive and $\det V_{n+1}>0$. Similarly, all the
principal minors of $U_{n+1}$ are those fo $U_n$ (hence positive), except for the last one:
\[ \det U_{n+1} = (v_{2(n+1)}-m_{2(n+1)}) \det U_n + C_2,
\]
so we have to demand that $v_{2(n+1)} > m_{2(n+1)} + \dfrac{C_2}{\det U_n}$ resulting in the positivity of the
matrix $U_{n+1}$ and $\det U_{n+1}>0$. Notice that $v_{2n+1}$ is arbitrary and we choose
$v_{2(n+1)} > \max (-\dfrac{C_1}{\det(V_n)},m_{2(n+1)} + \dfrac{C_2}{\det U_n})$, finishing the construction of
the sequences $v_n, u_n$. 

Finally, using Hamburger theorem, we find measures $\rho_+$ and $\rho_-$ such that
\[ v_n = \int\limits_{\mathbb R} \mu^n d\rho_+(\mu), \quad
u_n = \int\limits_{\mathbb R} \mu^n d\rho_-(\mu).
\]
As a result, for the signed measure $d\rho=d\rho_+ - d\rho_-$ it holds that
\[ \int\limits_{\mathbb R} \mu^n d\rho = \int\limits_{\mathbb R} \mu^n d\rho_+ -
\int\limits_{\mathbb R} \mu^n d\rho_- = v_n-u_n = m_n.
\]
\qed

Next theorem appears in \cite{bib:Boas1939}.
\begin{cor} \label{cor:PosSupp}Given a sequence of moments $H_n$ of the form \eqref{eq:HnSpcial}, there exists a function $S(\lambda)$,
realized in form \eqref{eq:SindrhoForm} with an analytic measure $d\bar\rho=\bbmatrix{d\rho_{11}&id\rho_{12}\\-id\rho_{12} & d\rho_{22}}$ 
as in Theorem \ref{thm:S0FromHn} with support on $[0,\infty]$.
\end{cor}
\Proof Using the construction of Theorem \ref{thm:S0FromHn} we can add a requirement on the choice of $v_{2n+1}, u_{2n+1}$ so that the conditions
of the Stieltjes moment problem are fulfilled. Additionally to the positivity of $V_n$ one has to require also that the matrices
\[ V'_n = \bbmatrix{
v_1&v_2&\ldots&v_{n+1}\\
v_2&v_2&\ldots&v_{n+2}\\
\vdots&\vdots&\ddots&\vdots\\
v_{n+1}&v_{n+2}&\ldots& v_{2n+1}}
\]
are positive (similar matrix $U'_n$ is constructed from $u_n$'s).
The condition is easily obtained when one considers the principal minors and uses the induction similarly to the proof of Theorem \ref{eq:HnSpcial}. \qed
\begin{thm} \label{thm:N0froS}
Suppose that a function $S(\lambda)$  possesses a realization \eqref{eq:SindrhoForm}
\[ S(\lambda) = I - \int\limits_{\mathbb R} \dfrac{1}{\lambda-i\mu} d\bar\rho(\mu)\sigma_1
\]
for an analytic measure $d\bar\rho = \bbmatrix{d\rho_{11}&id\rho_{12}\\-id\rho_{12} & d\rho_{22}}$. Then there exists an invertible node 
\[\mathfrak N_0 = \bbmatrix{
C_0 & A_\zeta, \mathbb X_0, A & B_0 & \sigma_1\\
\mathbb C^2 & \mathcal K & \mathbb C^2},
\]
where $A,A_\zeta$ are generators of analytic semi groups. The transfer function of the node $\mathfrak N_0$ is equal to $S(\lambda)$.
\end{thm}
\Proof We are going to explicitly construct such a node, based on \eqref{eq:SindrhoForm}.
Let us define two measures on column vector-functions with 2 entries:
\begin{equation} \label{eq:rho+rho-def}
 d\bar\rho^+ = \bbmatrix{d\rho_{11}^+ + d|\rho_{12}| & id\rho_{12}\\-id\rho_{12}&d\rho_{22}^+ + d|\rho_{12}|}, \quad
d\bar\rho^- = \bbmatrix{d\rho_{11}^- + d|\rho_{12}| &0\\0&d\rho_{22}^- + d|\rho_{12}|}.
\end{equation}
It is easy to see that $d\bar\rho = d\bar\rho^+ - d\bar\rho^-$ and we denote by $d|\bar\rho| = d\bar\rho^+ + d\bar\rho^-$.
Define a Hilbert space $\mathcal H$ of column vector-functions with 2 entries as follows
($\mathfrak R=supp(d\bar\rho)$)
\begin{equation} \label{eq:DefH4}
\mathcal H = \{ \bar u(\mu)=\bbmatrix{u_1(\mu)\\u_2(\mu)}\in\mathbb R^2 \mid \int_{\mathfrak R} \bar u^*(\mu) d|\bar\rho(\mu)| \bar u(\mu) < \infty \},
\end{equation}
equipped with the inner product
\begin{equation} \label{eq:DefH4innprod}
\langle \bar u, \bar v \rangle_{\mathcal H} =  \int_{\mathfrak R} \bar v^*(\mu) d|\bar\rho(\mu)| \bar u(\mu).
\end{equation}
This is a well defined object, because the measures $d\bar\rho^+, d\bar\rho^+$ are positive on the corresponding vector-functions.
The positivity of $d\bar\rho^-$ is immediate using \eqref{eq:rho+rho-def} on arbitrary integrable $\bbmatrix{v_1(\mu)\\v_2(\mu)}$:
\begin{multline*}
 \int_{\mathfrak R} \bbmatrix{v_1^*(\mu) & v_2^*(\mu)} \bbmatrix{d\rho_{11}^- + d|\rho_{12}| &0\\0&d\rho_{22}^- + d|\rho_{12}|} \bbmatrix{v(\mu) \\ v(\mu)} = \\
= \int_{\mathfrak R} |v_1|^2 (d\rho_{11}^- + d|\rho_{12}|) + \int_{\mathfrak R} |v_2|^2 (d\rho_{22}^- + d|\rho_{12}|) \geq 0.
\end{multline*}
For the measure $d\bar\rho^+$ we need more computations and the following formula can be shown:
\begin{multline*}
 \int_{\mathfrak R} \bbmatrix{v_1^*(\mu) & v_2^*(\mu)} \bbmatrix{d\rho_{11}^+ + d|\rho_{12}| & id\rho_{12}\\-id\rho_{12}&d\rho_{22}^+ + d|\rho_{12}|} \bbmatrix{v(\mu) \\ v(\mu)} = \\
= \int_{\mathfrak R} |v_1|^2 d\rho_{11}^+ +\int_{\mathfrak R} |v_2|^2 d\rho_{22}^+ + \int_{\mathfrak R} (\Re v_1+\Im v_2)^2 d\rho_{12}^+ + \int_{\mathfrak R}(\Re v_2-\Im v_1)^2 d\rho_{12}^+ \geq 0.
\end{multline*}
We define a Krein space $\mathcal K = \mathcal H$ as a set, equipped with the following sesqui-linear form
\begin{equation} \label{eq:DefK4innprod}
[\bar u,\bar v]_{\mathcal K} = \int_{\mathbb R} \bar v^*(\mu) d\bar\rho(\mu) \bar u(\mu).
\end{equation}

Define the operator $A=i\mu$ as the multiplication operator and 
\[ A_\zeta f = -i\mu f(\mu) - \sigma_1 \int_{\mathbb R} d\rho(\delta)f(\delta).\]
The operator $A_\zeta$ is a two-dimensional perturbation of the operator $-A$: each function $f(\mu)$ is mapped by $A_\zeta$
to the sum of $-i\mu f(\mu)$ and a constant function $K=-\sigma_1\int_{\mathbb R} d\rho(\delta)f(\delta)$.
The operators are generators of analytic semi-groups. Indeed, the group for $A$ is given by 
$e^{i\mu x}$ and is unitary. For the operator $A_\zeta$, we notice
that for big enough $\lambda>0$, we can explicitly write the inverse of $\lambda I - A_\zeta$. From
\[ (\lambda I - A_\zeta) f = (\lambda + i\mu)f(\mu) + \sigma_1 B_0^* f  = (\lambda + i\mu)f(\mu) +\sigma_1 K = g(\mu) \]
it follows that
\[ (\lambda I - A_\zeta)^{-1} g = \dfrac{g(\mu) - \sigma_1 K}{\lambda + i\mu},
\]
where the constant vector $K=K(\lambda,g)$ is found from the condition $K = B_0^* [\dfrac{g(\mu) - \sigma_1 K}{\lambda + i\mu}]$. Solving it we find that
\[ [I + \sigma_1 \int\limits_0^\infty \dfrac{d\bar\rho(\mu)}{\lambda+i\mu}] K = \int\limits_0^\infty \dfrac{d\bar\rho(\mu)g(\mu)}{\lambda+i\mu}
\]
and since for big enough $|\lambda|$ it holds that $\|\sigma_1 \int\limits_0^\infty \dfrac{d\bar\rho(\mu)}{\lambda+i\mu}\| < 1$, we obtain
\[ K = [I + \sigma_1 \int\limits_0^\infty \dfrac{d\bar\rho(\mu)}{\lambda+i\mu}]^{-1} \int\limits_0^\infty \dfrac{d\bar\rho(\mu)g(\mu)}{\lambda+i\mu}.
\]
From here it follows immediately that $(\lambda I - A_\zeta)^{-1}$ is bounded and $A_\zeta$ is a generator of an analytic semi group.

Obviously, $D(A_\zeta)=D(A)$. Define $\mathbb X_0 = I:\mathcal K\rightarrow\mathcal K$ - the identity operator. So, the conditions 
$\mathbb X_0(D(A)) = D(A) = \mathbb X_0^{-1}(D(A))$ of an invertible node are fulfilled.
Define $B_0 = C_0^* = I:\mathbb C^2\rightarrow\mathcal K$ and notice that $C_0:\mathcal K\rightarrow\mathbb C^2$ is an integration as follows ($f\in\mathcal K$)
\[ C_0 f u = \int_{\mathbb R} d\bar\rho(\mu) f(\mu).
\]
Then we compute for each $f\in D(A)$
\[ \begin{array}{lll}
A \mathbb X_0 f(\mu) + \mathbb X_0 A_\zeta f(\mu) + B_0\sigma_1C_0 f(\mu)=
i\mu f(\mu) - i\mu f(\mu)-\sigma_1 \int_{\mathbb R} d\rho(\delta)f(\delta) + \sigma_1 \int_{\mathbb R} d\bar\rho(\mu) \bar u(\mu) = 0,
\end{array} \]
which means that the set
\[\mathfrak N_0 = \bbmatrix{
I^* & A_\zeta, I, A & I & \sigma_1\\
\mathbb C^2 & \mathcal K & \mathbb C^2},
\]
is an invertible node. Its transfer function is
\[ \begin{array}{lll}
I - C_0 \mathbb X_0^{-1} (\lambda I - A)^{-1} B_0\sigma_1 
& = I - \int_{\mathbb R} I d\bar\rho(\mu) I (\lambda-i\mu)^{-1} I\sigma_1 \\
& = I - \int_{\mathbb R} \dfrac{1}{\lambda-i\mu} d\bar\rho(\mu)\sigma_1 \\
& = S(\lambda)
\end{array} \]
and the Theorem follows. \qed

\subsection{Construction of a vessel, realizing a given analytic potential}
Suppose that an analytic function $q(x)$ is given. We assume that $x_0 = 0$ for the simplicity of notations.
Using results of the previous Section \ref{sec:S0Constr}, or more precisely Corollary \ref{cor:PosSupp}, we construct an invertible node
$\mathfrak N_0$ to which we can apply the standard construction of a prevessel (see Section \ref{sec:preVConcstr}). We obtain in this manner a 
prevessel $\mathfrak{preV}$. Moreover, by Theorem \ref{thm:VFromS0} there exists an interval $\mathrm I$, including $x_0$ and a vessel
$\mathfrak V$ on $\mathrm I$, such that the potential of the vessel $q_V(x)$ exists and is analytic. Moreover, from the form of the zero moment
$H_0(x)$, by observing its $1,1$ entry we will obtain that $q_V(x) = q(x)$ on $\mathrm I$. So, if we are able to show that actually the vessel 
$\mathfrak V$ exists on the whole $\mathbb R$, we will realize the given potential by a vessel, constructing a scattering theory for it.

The following Theorem \ref{thm:SxrelizedSL} shows that there exists a transfer function $S(\lambda,x)$, which realizes the given potential.
This is a first sign
that a vessel $\mathfrak V$ realizing $q(x)$ on $\mathbb R$ exists. Let us denote by $\Phi(\mu,x), \Phi_*(\mu,x)$ the fundamental solutions
of \eqref{eq:InCC}, \eqref{eq:OutCC} respectively. First we notice that the columns and rows of 
the fundamental matrices are in $\mathcal K$. Indeed
\begin{equation}\label{eq:PhiForm}
 \Phi(i\mu,x) = \bbmatrix{
\cos(\sqrt{\mu}x) & -i\sqrt{\mu} \sin(\sqrt{\mu} x)\\
-i\dfrac{\sin(\sqrt{\mu} x)}{\sqrt{\mu}} & \cos(\sqrt{\mu}x)}
\end{equation}
and it is obvious that $| \Phi(i\mu,x)| < \sqrt{\mu} C$ for some constant. Since $d|\rho(\mu)|$ is analytic the integral
\[ \int_0^\infty |\Phi(i\mu,x)|^2 d|\rho(\mu)| < C^2 \int_0^\infty \mu d|\rho(\mu)|
\]
is finite, which means that the columns and rows of $\Phi(i\mu,x)$ are in $\mathcal K$. To prove that the columns and rows of $\Phi_*(i\mu,x)$ are in $\mathcal K$ for each
$x\in\mathbb R$, we need to learn its structure first. From \eqref{eq:ysimpleODEs}, \eqref{eq:qbeta} it follows that
\[ \Phi_*(i\mu,x) = \bbmatrix{\phi(\mu,x)& i \psi(\mu,s)\\
-i(\dfrac{\partial}{\partial x} \phi(\mu,x)-\beta(x)\phi(\mu,x)) & \dfrac{\partial}{\partial x} \psi(\mu,x)-\beta(x)\psi(\mu,x)}
\]
where $\beta(x)=-\dfrac{1}{2} \int\limits_0^x q(y)dy$ and $\phi(\mu,x), \psi(\mu,x)$ are solutions of \eqref{eq:SL} with the initial conditions
\[ \phi(\mu,0) = 1, \quad \dfrac{\partial}{\partial x} \phi(\mu,0)=0,\quad
\psi(\mu,0) = 0, \quad \dfrac{\partial}{\partial x} \psi(\mu,0)=1.
\]
The structure of these solutions is very well known \cite{bib:FaddeyevII}. Using variation of coefficients they satisfy
\[ \begin{array}{ll}
\phi(\mu,x) = \cos(\sqrt{\mu}x) + \int\limits_0^x \dfrac{\sin(\sqrt{\mu}(x-y))}{\sqrt{\mu}} q(y) \phi(\mu,y)dy, \\
\psi(\mu,x) = \sin(\sqrt{\mu}x) + \int\limits_0^x \dfrac{\sin(\sqrt{\mu}(x-y))}{\sqrt{\mu}} q(y) \psi(\mu,y)dy.
\end{array} \]
And from their Liouville–Neumann series solutions we obtain
\begin{multline*}
\phi(\mu,x) = \cos(\sqrt{\mu}x) - \int\limits_0^x \dfrac{\sin(\sqrt{\mu}(x-y))}{\sqrt{\mu}} q(y) \cos(\sqrt{\mu}y)dy + \\
+ \int\limits_0^x \dfrac{\sin(\sqrt{\mu}(x-y))}{\sqrt{\mu}}  q(y) \int\limits_0^y \dfrac{\sin(\sqrt{\mu}(y-y_1))}{\sqrt{\mu}}  q(y_1) \cos(\sqrt{\mu}y)dy_1  dy - \cdots
\end{multline*}
Since $|\dfrac{\sin(\sqrt{\mu}(x-y))}{\sqrt{\mu}}|, |\cos(\sqrt{\mu}y)|< 1$, $|q(y)| < M_x$ on $[0,x]$ (or on $[x,0]$ for $x<0$) for a constant $M_x>0$
\[ \sup_\mathbb R|\phi(\mu,x)| < 1 + \int_0^{|x|} M_x dy + \int\limits_0^{|x|} \int\limits_0^y M_x^2 dy_1 dy + \cdots = e^{M_x |x|} < \infty.
\]
Differentiating the formula for $\phi(\mu,x)$ we find that
\[ \dfrac{\partial}{\partial x} \phi(\mu,x) = -\sqrt{\mu}\sin(\sqrt{\mu}x) + \int\limits_0^x \cos(\sqrt{\mu}(x-y)) q(y) \phi(\mu,y)dy,
\]
from where it follows that $|\dfrac{\partial}{\partial x} \phi(\mu,x)| < C \sqrt{\mu}$, knowing the bound for $|\phi(\mu,x)|$. Similarly,
one finds that $\phi(\mu,x)$ and its $x$-derivative satisfy the same bounds and as a result $|\Phi_*(\mu,x)| < \sqrt{\mu}C$ and the columns and rows of $\Phi_*(\mu,x)$ are in $\mathcal K$.
\begin{thm}[Transfer function construction] \label{thm:SxrelizedSL}Let $q(x)$ be an analytic functions and let $H_n(x)$ be the moments, constructed in \eqref{eq:Hi+1HiFs}.
Suppose that $S(\lambda)$ is realized in the form \eqref{eq:SindrhoForm}
\[ S(\lambda) = I - \int\limits_{\mathbb R} \dfrac{1}{\lambda-i\mu} d\bar\rho(\mu)\sigma_1,
\]
with an analytic measure $d\bar\rho$ and satisfies $H^S_n=H_n(0)$. Let $\Phi(\lambda,x),\Phi_*(\lambda,x)$ be the fundamental solutions of 
\eqref{eq:InCC}, \eqref{eq:OutCC} respectively. Then the function
\[ S(\lambda,x) = I - \int\limits_{\mathbb R} \dfrac{1}{\lambda-i\mu} \Phi_*(i\mu,x) d\bar\rho(\mu) \Phi^*(i\mu,x) \sigma_1,
\]
satisfies $H_n^S(x)=H_n(x)$ for all $x\in\mathbb R$.
\end{thm}
\Proof Differentiating the $n$-th moment of $S(\lambda,x)$ 
\[ H^S_n(x) =  \int\limits_{\mathbb R}  \Phi_*(i\mu,x) d\bar\rho(\mu) (i\mu)^n\Phi^*(i\mu,x)
\]
we find that
\begin{multline*}
\dfrac{\partial}{\partial x} H^S_n(x) = \int\limits_{\mathbb R}  \sigma_1^{-1}(\sigma_2i\mu+\gamma_*(x))\Phi_*(i\mu,x) d\bar\rho(\mu) (i\mu)^n\Phi^*(i\mu,x) - \\
- \int\limits_{\mathbb R}  \Phi_*(i\mu,x) d\bar\rho(\mu) (i\mu)^n\Phi^*(i\mu,x)  (\sigma_2i\mu+\gamma)\sigma_1^{-1} = \\
= \sigma_1^{-1}\sigma_2H^S_{n+1} - H^S_{n+1} \sigma_2\sigma_1^{-1} + \sigma_1^{-1} \gamma_* H^S_n - H^S_n \gamma\sigma_1^{-1}
\end{multline*}
which is identical to \eqref{eq:DHn}. So $H_n^S(x)$ and $H_n(x)$ have the same initial conditions and satisfy the same differential equations, so they are identical
by the uniqueness of the moments Theorem \ref{thm:momentUnique}.
\qed
\begin{cor} \label{cor:HnfromPhiPhi*} $H_n(x) =  \int\limits_{\mathbb R}  \Phi_*(i\mu,x) d\bar\rho(\mu) (i\mu)^n\Phi^*(i\mu,x)$. Particularly, for $n=0$
\[ H_0(x) = \int\limits_{\mathbb R}  \Phi_*(i\mu,x) d\bar\rho(\mu) \Phi^*(i\mu,x).
\]
\end{cor}
Assume that for a given analytic $q(x)$ we have constructed moments $H_n(x)$, a measure $d\bar\rho$ in Corollary \ref{cor:PosSupp} and a node
$\mathfrak N_0$ (Theorem \ref{thm:N0froS}). Applying the standard construction Theorem \ref{thm:STConpreVessel} to $\mathfrak N_0$, we obtain a prevessel:
\[ \begin{array}{lllllll}
\mathfrak{preV} & = \bbmatrix{C(x) & A_\zeta, \mathbb X(x), A & B(x)&  \sigma_1,\sigma_2,\gamma \\
\mathbb C^2 & \mathcal{K} & \mathbb C^2}.
\end{array} \]
In order to show that the operator $\mathbb X(x)$ is globally defined, we construct an ``inverse vessel'' as follows. We use similar to \eqref{eq:StConB}, \eqref{eq:StConC} 
definitions, using the fundamental matrix $\Phi_*(\lambda,x)$ instead of $\Phi(\lambda,x)$:
\begin{defn} Define the operators
\begin{eqnarray}
B_*(x) = \dfrac{1}{2\pi i} \int\limits_\Gamma (\lambda I + A_\zeta)^{-1} B_0 \Phi^*_*(-\bar\lambda,x) d\lambda,, \\
C_*(x) = \dfrac{1}{2\pi i} \int\limits_\Gamma \Phi_*(\lambda,x) C_0 (\lambda I - A)^{-1} d\lambda, \\
\label{eq:DefX*} \mathbb X_*(x) = I - \int_0^x B_*(x)\sigma_2C_*(y)dy.
\end{eqnarray}
\end{defn}
\begin{lemma} The operators $B_*(x), C_*(x), \mathbb X_*(x)$ satisfy
\begin{eqnarray}
\label{eq:DB*} \dfrac{d}{dx} B_*(x) = (A_\zeta B_*(x)\sigma_2 - B_*(x)\gamma_*(x))\sigma_1^{-1}, \quad B_*(0) = B_0, \\
\label{eq:DC*} \dfrac{d}{dx} C_*(x)u = \sigma_1^{-1}(\sigma_2 C(x) A + \gamma_*(x) C_*(x)) u, \quad u\in D(A), \quad C_*(0) = C_0, \\
\label{eq:X*Lyapunov} A_\zeta \mathbb X_*(x) u + \mathbb X_*(x)A u + B_*(x)\sigma_1C_*(x) u = 0,\quad u\in D(A).
\end{eqnarray}
\end{lemma}
\Proof Immediate from the definitions. The Lyapunov equation \eqref{eq:X*Lyapunov} follows similarly to the proof of Lemma \ref{lemma:Redund}.
\begin{lemma} Define moments $G_n(x) = C_*(x)A^n B(x)$. Then $G_n(x) = H_n(x)$, particularly 
\begin{equation} \label{lemma:C8BisH0} 
C_*(x)B(x)=H_0(x).
\end{equation}
\end{lemma}
\Proof The moments $G_n(x)$ are well defined, since $A^n B(x)$ is an element of $\mathcal K$ for each $n$. Then
\begin{multline*}
 \dfrac{d}{dx} G_n(x) = \dfrac{d}{dx} [C_*(x)A^n B(x)] = \\
 = \sigma_1^{-1} (\sigma_2 C_*(x)A^{n+1} + \gamma_*(x) C_*(x)A^n) B(x) - C_*(x) A^n (A B(x)\sigma_2 + B(x)\gamma)\sigma_1^{-1} = \\
= \sigma_1^{-1} (\sigma_2 G_{n+1}(x) + \gamma_*(x) G_n(x)) - (G_{n+1}(x) \sigma_2 + G_n(x)\gamma)^{-1} = \\
=  	\sigma_1^{-1}\sigma_2G_{n+1}(x) - G_{n+1}(x) \sigma_2\sigma_1^{-1} + \sigma_1^{-1} \gamma_*(x) G_n(x) - G_n(x) \gamma\sigma_1^{-1},
\end{multline*}
which coincides with \eqref{eq:DHn}. Moreover, $G_n(0) = H_n(0) = C_0A^n B_0$ by their constructions. So, by the uniqueness of the moments Theorem
\ref{thm:momentUnique}, $G_n(x)=H_n(x)$. \qed

From the equation \eqref{lemma:C8BisH0} we obtain relations between the operators:
\begin{thm} \label{thm:vesselmappings}
The following formulas hold
\begin{eqnarray*}
\mathbb X(x) B_*(x) = B(x),  \quad C_*(x)\mathbb X(x) = C(x), \\
\mathbb X_*(x) B(x) = B_*(x),  \quad C(x)\mathbb X_*(x) = C_*(x).
\end{eqnarray*}
\end{thm}
\Proof Let us prove the identity $C_*(x)\mathbb X(x) = C(x)$ and the rest are obtained in a similar manner. The following identities are applied to
an element $u\in D(A)$:
\[ \begin{array}{lllllll} 
\dfrac{d}{dx} [C_*(x)\mathbb X(x)] & =  \sigma_1^{-1} [ \sigma_2 C_*(x) A + \gamma_*(x)C_*(x)] \mathbb X(x) + C_*(x) B(x)\sigma_2 C(x)  \\
& =  \sigma_1^{-1} \sigma_2 C_*(x) A \mathbb X(x) + \gamma_*(x) C_*(x) \mathbb X(x) + H_0(x) \sigma_2 C(x)\\
& = \text{ using \eqref{eq:Lyapunov}} \\
& = \sigma_1^{-1} \sigma_2 C_*(x) [-\mathbb X(x)A_\zeta - B(x) \sigma_1C(x)] + \sigma_1^{-1} \gamma_*(x) C_*(x) \mathbb X(x) + H_0(x) \sigma_2 C(x)\\
& = - \sigma_1^{-1} \sigma_2 C_*(x) \mathbb X(x)A_\zeta - \sigma_1^{-1} \sigma_2 C_*(x) B(x) \sigma_1C(x) + \\
& \quad \quad + \sigma_1^{-1} \gamma_*(x) C_*(x) \mathbb X(x) + H_0(x) \sigma_2 C(x)\\
& = \text{using \eqref{lemma:C8BisH0} and \eqref{eq:Linkage} } \\
& = - \sigma_1^{-1} \sigma_2 [C_*(x) \mathbb X(x)] A_\zeta + \sigma_1^{-1} \gamma_*(x) [C_*(x) \mathbb X(x)] - [\gamma_*(x)-\gamma] C(x)
\end{array} \]
in other words the operator $C_*(x)\mathbb X(x)$ satisfies the following non-homogeneous differential equation
\[ \dfrac{d}{dx} Y = - \sigma_1^{-1} \sigma_2 Y A_\zeta + \sigma_1^{-1} \gamma_*(x) Y - [\gamma_*(x)-\gamma] C(x).
\]
On the other hand, $C(x)$ satisfies the same differential equation:
\[ \dfrac{d}{dx} C(x) = - \sigma_1^{-1} \sigma_2 C(x) A_\zeta + \sigma_1^{-1} \gamma_*(x) C(x) - [\gamma_*(x)-\gamma] C(x) = \sigma_1^{-1}(- \sigma_2 C(x) A_\zeta + \gamma C(x)).
\]
Since $C_*(0)\mathbb X(0) = C(0) = C_0$ the result follows, by the uniqueness of the solution. Similarly, using \eqref{eq:DB}, \eqref{eq:DefX*} and \eqref{eq:X*Lyapunov}
\[ \dfrac{d}{dx} [\mathbb X_*(x)B(x)] = A_\zeta [\mathbb X_*(x)B(x)] \sigma_2\sigma_1^{-1} - [\mathbb X_*(x)B(x)] \gamma\sigma_1^{-1} - B_*(x)[\gamma_*(x)-\gamma]\sigma_1^{-1},
\]
and $\mathbb X_*(x)B(x)$, substituted with $B_*(x)$ satisfies the same differential equation. Thus $\mathbb X_*(x)B(x)=B_*(x)$. As a result, we obtain that
\[ C(x)B_*(x) = C(x)\mathbb X_*(x)B(x) = C_*(x)B(x) = H_0(x).
\]
Then the equations $C(x)\mathbb X_*(x) = C_*(x)$,  $\mathbb X(x) B_*(x) = B(x)$ follows in the same manner. \qed

\begin{cor} \label{cor:XxInvertible}
The operator $\mathbb X(x)$ is invertible for all $x\in\mathbb R$ with the inverse $\mathbb X_*(x)$.
\end{cor}
\Proof Notice that from Theorem \ref{thm:vesselmappings} it follows that
\begin{eqnarray*} \dfrac{d}{dx} [\mathbb X(x)\mathbb X_*(x)] = B(x)\sigma_2 C(x) \mathbb X_*(x) - \mathbb X(x) B_*(x)\sigma_2 C_*(x) = 
B(x)\sigma_2 C_*(x) - B(x)\sigma_2 C_*(x) = 0, \\
\dfrac{d}{dx} [\mathbb X_*(x)\mathbb X(x)] = -B_*(x)\sigma_2 C_*(x) \mathbb X(x) - \mathbb X_*(x) B(x)\sigma_2 C(x) =
B_*(x)\sigma_2 C(x) - B_*(x)\sigma_2 C(x) = 0.
\end{eqnarray*}
Since $\mathbb X(0)\mathbb X_*(0) = \mathbb X_*(0)\mathbb X(0) = I$, the Corollary follows. 
\qed
\begin{thm} \label{thm:VesselFromqOnR}
The vessel,
\[ \begin{array}{lllllll}
\mathfrak{V} & = \bbmatrix{C(x) & A_\zeta, \mathbb X(x), A & B(x)&  \sigma_1,\sigma_2,\gamma, \gamma_*(x) \\
\mathbb C^2 & \mathcal{K} & \mathbb C^2&\mathbb R}
\end{array} \]
obtained by applying the standard construction to the node $\mathfrak N_0$ exists on $\mathbb R$.
The matrix function $\gamma_*(x)$, defined by the linkage condition \eqref{eq:Linkage} realizes the potential $q(x)$ on $\mathbb R$.
\end{thm}
\Proof Since $\mathbb X(x)$ is globally invertible, we can repeat the proof of Theorem \ref{thm:VFromS0} in order to show that
the prevessel $\mathfrak{preV}$ is an invertible node. this shows that $\mathfrak{V}$ is a vessel, realizing an analytic potential
$q_V(x)$ on $\mathbb R$, since $\mathbb X(x)$ is globally invertible. Then by theorem \ref{thm:SxrelizedSL} the moments of the vessel
$\mathfrak{V}$ are equal to the moments $H_n(x)$, particularly $H_0(x) = C(x)\mathbb X^{-1}(x)B(x)$ for which the $1,1$ entry means
that $q_V(x) = q(x)$.
\qed

\section{KdV evolutionary vessels}
Let us evolve a SL vessel with respect to $t$. Some of the results presented here can be found in
\cite{bib:KdVVessels, bib:KdVHierarchy} for symmetric vessels.

We consider the following notion
\begin{defn} The collection of operators and spaces 
\begin{equation} \label{eq:DefKdVpreV}
\mathfrak{preV}_{KdV} = \bbmatrix{C(x,t) & A_\zeta, \mathbb X(x,t), A & B(x,t)&  \sigma_1,\sigma_2,\gamma \\
\mathbb C^2 & \mathcal{K} & \mathbb C^2}
\end{equation}
is called a \textbf{KdV preVessel}, if the following conditions hold:
1. $\mathfrak{preV}_{KdV}$ is a node for all $x,t\in\mathbb R$, 2.
operator $B(x,t)\sigma_2$ is $A^2$-regular, $B(x,t)\gamma$ is $A$-regular 3. $C(x,t), \mathbb X(x,t), B(x,t)$ are differentiable in
both variables, when the other one is fixed, subject to the conditions \eqref{eq:DB}, \eqref{eq:DC}, \eqref{eq:DX} and the following evolutionary equations
(for arbitrary $u\in D(A), v\in D(A)$)
\begin{eqnarray}
\label{eq:DBt} \frac{\partial}{\partial t} B  & = iA \dfrac{\partial}{\partial x} B & = - i A (A B \sigma_2 + B \gamma) \sigma_1^{-1}, \\
\label{eq:DCt} \frac{\partial}{\partial t} C u & = - i \dfrac{\partial}{\partial x} C A_\zeta u & 
= - i \sigma_1^{-1}  (-\sigma_2 C A_\zeta + \gamma C) A_\zeta u, \\
\label{eq:DXt} \frac{\partial}{\partial t} \mathbb X v & =  i (A \dfrac{\partial}{\partial x}\mathbb X - i \dfrac{\partial}{\partial x}\mathbb X A_\zeta+iB\gamma C)v &
= i (A B \sigma_2 C - i B\sigma_2 C A_\zeta+iB\gamma C)v ,
\end{eqnarray}
where $\sigma_2=\sigma_2^*$, $\gamma^*=-\gamma$ are $2\times 2$ matrices. The prevessel $\mathfrak{preV}$ is called symmetric if
$A_\zeta=A^*$ and $C(x,t)=B^*(x,t)$ for all $x,t\in\mathbb R$.
\end{defn}

\begin{defn} \label{def:KdVVessel} The collection of operators, spaces and an open set $\Omega\subseteq\mathbb R^2$ 
\begin{equation} \label{eq:DefKdVVGen}
\mathfrak{V}_{KdV} = \bbmatrix{C(x,t) & A_\zeta, \mathbb X(x,t), A & B(x,t)&  \sigma_1,\sigma_2,\gamma,\gamma_*(x,t) \\
\mathbb C^2 & \mathcal{K} & \mathbb C^2 & \Omega }
\end{equation}
is called a (non-symmetric) \textbf{KdV vessel}, if $\mathfrak{V}_{KdV}$ is a KdV prevessel, $\mathbb X(x,t)$ is invertible on $\Omega$, $\mathfrak{V}_{KdV}$
is also an invertible node.
The $2\times 2$ matrix-function $\gamma_*(x,t)$ satisfies the linkage condition \eqref{eq:Linkage}.
The vessel $\mathfrak{V}_{KdV}$ is called symmetric if $A_\zeta=A^*$ and $C(x,t)=B^*(x,t)$ for all $x,t\in\Omega$.
\end{defn}
\begin{thm} Let $\mathfrak{V}_{KdV}$ be a KdV vessel. Suppose that the moments $H_0,\ldots,H_{n+1}$ are finite and
differentiable, then
\begin{equation}\label{eq:DHntKdV}
\dfrac{\partial}{\partial t} H_n = i \dfrac{\partial}{\partial x} H_{n+1} + i \dfrac{\partial}{\partial x} [H_0] \sigma_1 H_n.
\end{equation}
The transfer function $S(\lambda,x,t)$ \eqref{eq:S0realized} satisfies the following differential equation
\begin{equation} \label{eq:PDEforS}
 \dfrac{\partial}{\partial t} S(\lambda,x,t) = i\lambda \dfrac{\partial}{\partial x} S(\lambda,x,t) +
i \dfrac{\partial}{\partial x} [H_0] \sigma_1 S(\lambda,x,t).
\end{equation}
\end{thm}
\Proof Consider the formula for the moments first.
\[ \begin{array}{lllllll}
\dfrac{\partial}{\partial t} H_n & = \dfrac{\partial}{\partial t}[C\mathbb X^{-1}A^n B] = C_t\mathbb X^{-1}A^nB - 
C \mathbb X^{-1} \mathbb X_t \mathbb X^{-1}A^nB + C\mathbb X^{-1}A^n B_t = \\
& = \text{using evolutionary conditions \eqref{eq:DBt}, \eqref{eq:DXt} } \\
& = C_x(-iA_\zeta) \mathbb X^{-1}A^nB - C \mathbb X^{-1} (iA \mathbb X_x-i\mathbb X_x A_\zeta + iB\gamma C) \mathbb X^{-1}A^nB + C\mathbb X^{-1}A^n (iA) B_x = \\
& = \text{using \eqref{eq:DB}, \eqref{eq:Linkage} and \eqref{eq:Lyapunov}} \\
& = i \dfrac{\partial}{\partial x} H_{n+1} + i \dfrac{\partial}{\partial x} [H_0] \sigma_1 H_n,
\end{array} \]
Similarly one shows the formula \eqref{eq:PDEforS}.
\qed
\begin{cor} The potential $\gamma_*(x,t)$ of a KdV vessel satisfies the following differential equation
\begin{equation} \label{eq:Dgamma*tKdV}
(\gamma_*)_t = - i \gamma_* (H_0)_x\sigma_1 + i \sigma_1 (H_0)_{xx} \sigma_1 +i \sigma_1 (H_0)_x \gamma_*.
\end{equation} 
\end{cor}
\Proof From the linkage condition and \eqref{eq:DHntKdV} for $n=0$ it follows that
\[ \begin{array}{lllllll}
(\gamma_*)_t & = \sigma_2 (H_0)_t \sigma_1 - \sigma_1 (H_0)_t \sigma_2 = \\
& = \sigma_2 [i(H_1)_x  + i (H_0)_x \sigma_1 H_0] \sigma_1 - \sigma_1 [i (H_1)_x + i (H_0)_x \sigma_1 H_0] \sigma_2 \\
& = i \sigma_1 [\sigma_1^{-1} \sigma_2 (H_1)_x - (H_1)_x \sigma_2\sigma_1^{-1}]\sigma_1 + 
i \sigma_2 (H_0)_x \sigma_1 H_0 \sigma_1 - i \sigma_1 (H_0)_x \sigma_1 H_0 \sigma_2.
\end{array} \]
For the first term in this expression we can use the formula \eqref{eq:DHn} for $n=0$, 
then
\[ \begin{array}{lllllll}
(\gamma_*)_t & = i \sigma_1 \dfrac{\partial}{\partial x}[(H_0)_x - \sigma_1^{-1}\gamma_*H_0 + H_0 \gamma \sigma_1^{-1} ]\sigma_1 + 
i \sigma_2 (H_0)_x \sigma_1 H_0 \sigma_1 - i \sigma_1 (H_0)_x \sigma_1 H_0 \sigma_2 \\
& = i \sigma_1(H_0)_{xx}\sigma_1 - i \dfrac{\partial}{\partial x} [\gamma_*H_0 \sigma_1 + \sigma_1 H_0 \gamma]+ 
i \sigma_2 (H_0)_x \sigma_1 H_0 \sigma_1 - i \sigma_1 (H_0)_x \sigma_1 H_0 \sigma_2 \\
& = i \sigma_1(H_0)_{xx}\sigma_1 - i \gamma_* (H_0)_x\sigma_1 - i (\gamma_*)_x H_0 \sigma_1 + \sigma_1 (H_0)_x \gamma +  i \sigma_2 (H_0)_x \sigma_1 H_0 \sigma_1 - i \sigma_1 (H_0)_x \sigma_1 H_0 \sigma_2.
\end{array} \]
Then notice that 
\[ \begin{array}{lllllll}
- i (\gamma_*)_x H_0 \sigma_1 + \sigma_1 (H_0)_x \gamma +  i \sigma_2 (H_0)_x \sigma_1 H_0 \sigma_1 - i \sigma_1 (H_0)_x \sigma_1 H_0 \sigma_2 = \\
= - i [\sigma_2 (H_0)_x\sigma_1 - \sigma_1 (H_0)_x\sigma_2] H_0 \sigma_1 + \sigma_1 (H_0)_x \gamma +  i \sigma_2 (H_0)_x \sigma_1 H_0 \sigma_1 - i \sigma_1 (H_0)_x \sigma_1 H_0 \sigma_2 = \\
= i \sigma_1 (H_0)_x [\gamma + \sigma_2 H_0 \sigma_1 - \sigma_1 H_0 \sigma_2] \\
= i \sigma_1 (H_0)_x \gamma_*,
\end{array} \]
and the result follows. \qed
\begin{cor} \label{cor:H0toKdV}The potential $q(x,t)$ of a KdV vessel satisfies the KdV equation \eqref{eq:KdV} on $\Omega$, which is equivalent to \eqref{eq:Dgamma*tKdV}.
\end{cor}
\Proof From the linkage condition it follows that for $\beta = \dfrac{1}{2} \int_0^x q(y,t)dy$
\[ \gamma_* = \bbmatrix{-i(\beta'_x-\beta^2)&-\beta\\\beta&i}.
\]
Moreover, by \eqref{eq:qbeta}, the KdV equation for $q(x,t)$ follows from the
differential equation for $\beta(t,x)$:
\begin{equation} \label{eq:PDEbeta}
 - 4 \beta_t = -6 (\beta_x)^2 + \beta_{xxx}.
\end{equation}
Then using \eqref{eq:Dgamma*tKdV} and \eqref{eq:Linkage}
\[ \begin{array}{lllllll}
-\beta_t & = \bbmatrix{1&0} (\gamma_*)_t \bbmatrix{0\\1} =
\bbmatrix{1&0} [- i \gamma_* (H_0)_x\sigma_1 + i \sigma_1 (H_0)_{xx} \sigma_1 +i \sigma_1 (H_0)_x \gamma_*] \bbmatrix{0\\1} = \\
& = \bbmatrix{1&0} [- i (\gamma+\sigma_2H_0\sigma_1-\sigma_1H_0\sigma_2)(H_0)_x\sigma_1 + i \sigma_1 (H_0)_{xx} \sigma_1 +i \sigma_1 (H_0)_x \gamma_*] \bbmatrix{0\\1}.
\end{array} \]
We have seen in \eqref{eq:H0Mat} that
\[ H_0 = \bbmatrix{ -\beta & \dfrac{-i\pi_{11}}{2} \\ \dfrac{i\pi_{11}}{2} & h_0^{22}}. \]
Plugging this formula into the last expression we will find that
\[ -\beta_t = -(h_0^{22})_x-\beta_x\beta^2 + 2 (\beta_x)^2 + \beta\beta_{xx} - \dfrac{1}{2} \beta_{xxx}.
\]
From the last formula of \eqref{eq:Hi+1HiFs} it follows that
\[ (h_0^{22})_x = \dfrac{1}{2i} [(\dfrac{-i\pi_{11}}{2})_{xx}-2\beta(\dfrac{-i\pi_{11}}{2})_x ]=
-\beta^2\beta_x + \dfrac{1}{2} (\beta_x)^2 + \beta\beta_{xx} - \dfrac{1}{4} \beta_{xxx}.
\]
Plugging this expression into the formula for $-\beta_t$, we will obtain \eqref{eq:PDEbeta}
\[ -\beta_t = \dfrac{3}{2} (\beta_x)^2 - \dfrac{1}{4} \beta_{xxx},
\]
which the KdV equation for $\beta(x)$. Differentiating it with respect to $x$, we will obtain the regular KdV equation \eqref{eq:KdV}
for $q(x,t)=2\beta_x(x,t)$. It is a matter of simple algebraic calculations to verify that the $1,1$ entry of \eqref{eq:Dgamma*tKdV} is equivalent to
\eqref{eq:KdV}, since the $1,1$ entry of $\gamma_*(x,t) = -i(\beta_x-\beta^2)$ is expressible in terms of $\beta$. \qed

Now we obtain the Main Theorem, because the fact that $\mathbb X(x)$ is invertible for a fixed $x$ implies that its norm is bounded from below and
a small perturbation of it is still invertible.
\begin{mthm} \label{mthm}Suppose that $q(x)$ is an analytic function on $\mathbb R$. There exists a KdV vessel, which exists on $\Omega \subseteq \mathbb R^2$. 
For each $x\in\mathbb R$ there exists $T_x>0$ such that $\{x\}\times[-T_x,T_x]\in\Omega$. The potential $q(x)$ is realized by the vessel for $t=0$. 
\end{mthm}
\Proof For the given analytic potential $q(x)$ we construct a SL vessel
\[ \begin{array}{lllllll}
\mathfrak{V} & = \bbmatrix{C(x) & A_\zeta, \mathbb X(x), A & B(x)&  \sigma_1,\sigma_2,\gamma, \gamma_*(x) \\
\mathbb C^2 & \mathcal{K} & \mathbb C^2} \\
& B(x) = \Phi^*(i\mu,x), \quad C(x) = \dfrac{1}{2\pi i} \int\limits_\Gamma \Phi(\lambda,x) C_0 (\lambda I + A_\zeta)^{-1} d\lambda \\
& \mathbb X(x) = I + \int_0^x B(y)\sigma_2C(y)dy = \\
&\quad = I + \int_\Gamma \dfrac{\Phi^*(i\mu,x)\sigma_1\Phi(\lambda,x)-\sigma_1}{(\lambda-i\mu)}C_0(\lambda I + A_\zeta)^{-1}d\lambda ,
\end{array} \]
defined in Theorem \ref{thm:VesselFromqOnR}. The last formula for $\mathbb X(x)$ comes from an easily checkable fact that
$\dfrac{\partial}{\partial x} \dfrac{\Phi^*(i\mu,x)\sigma_1\Phi(\lambda,x)-\sigma_1}{(\lambda-i\mu)} = \Phi^*(i\mu,x)\sigma_2\Phi(\lambda,x)$.
On its basis we define A KdV vessel \eqref{eq:DefKdVVGen}
\[ \mathfrak{V}_{KdV} = \bbmatrix{C(x,t) & A_\zeta, \mathbb X(x,t), A & B(x,t)&  \sigma_1,\sigma_2,\gamma,\gamma_*(x,t) \\
\mathbb C^2 & \mathcal{K} & \mathbb C^2 & \Omega} \]
as follows:
\[ \begin{array}{lllllll}
B(x,t) = \Phi^*(i\mu,x - \mu t), \quad C(x,t) = \dfrac{1}{2\pi i} \int\limits_\Gamma \Phi(\lambda,x-i\lambda t) C_0 (\lambda I + A_\zeta)^{-1} d\lambda, \\
\mathbb X(x,t) =  I + \int_\Gamma \dfrac{\Phi^*(i\mu,x-\mu t)\sigma_1\Phi(\lambda,x-i\lambda t)-\sigma_1}{(\lambda-i\mu)}C_0(\lambda I + A_\zeta)^{-1}d\lambda.
\end{array} \]
It is a matter of simple algebraic calculations to verify that $B(x,t), C(x,t), \mathbb X(x,t)$ satisfy the conditions of a KdV prevessel.
One has to use the fact that
\begin{multline*}
 \dfrac{\partial}{\partial t} \dfrac{\Phi^*(i\mu,x-\mu t)\sigma_1\Phi(\lambda,x-i\lambda t)-\sigma_1}{(\lambda-i\mu)} = \\
 = i \Phi^*(i\mu,x-\mu t)\sigma_2\Phi(\lambda,x-i\lambda t) (\lambda+i\mu) + i \Phi^*(i\mu,x-\mu t)\gamma \Phi(\lambda,x-i\lambda t).
\end{multline*}

Finally, notice that for $t=0$ the operator $\mathbb X(x,0)$ equals to the operator $\mathbb X(x)$ of the SL vessel, constructed for $q(x)$. Thus the set $\Omega$, on which
$\mathbb X(x,t)$ is invertible includes $\mathbb R\times \{0\}$. Moreover, since $(A_\zeta^*)^n C^*(x)\sigma_2$ is 
exists for all $n$ by the construction, we obtain that $\sigma_2C(x)A_\zeta$ is a well defined BOUNDED functional on $\mathcal K$. Thus the expression
\[ A B(x,s) \sigma_2 C(x,s) - i B(x,s)\bbmatrix{1\\0}[A^*_\zeta C^*(x) \bbmatrix{1\\0}]^* +iB(x,s)\gamma C(x,s) \]
is a bounded operator on $\mathcal K$. As a result, the operator $\mathbb X(x,t)$ is bounded for some $t\in[-T_x,T_x]$, where
\[ 0<T_x < \dfrac{\| A B(x,s) \sigma_2 C(x,s) - i B(x,s)\sigma_2 C(x,s) A_\zeta+iB(x,s)\gamma C(x,s)\|}{\|\mathbb X^{-1}(x,0)\|}. \qed \]
Finally, we present a Theorem, providing a conclusion that this theory of vessels is the ultimate tool for studying solutions of \eqref{eq:KdV}.
A most general Theorem in this connection is to show that if there is an open set $\Omega$, where solution is known to exist, then the operator $\mathbb X(x,t)$ is invertible
in this region. Since we use the uniqueness of solutions for ODEs, the set $\Omega$ must be at least simply-connected. On the other hand,
we do not want enter into topological difficulties, arising from such a general assumption, so we choose a very important and practical
case of a strip. So, if it is known that the solution of \eqref{eq:KdV} exists on $\mathbb R\times [0,T]$, we would like to show that the vessel,
which realizes $q(x)$ at $t=0$ will exist on this strip.

The idea of the proof of such a Theorem is very simple. We actually explicitly construct the inevsrse of $\mathbb X(x,t)$, using the assumption
that that there exists a solution $q(x,t)$ of \eqref{eq:KdV} on $\mathbb R\times [0,T]$. For a potential $q(x,t)$ on the strip, which solves \eqref{eq:KdV}
we can define 
\[ \beta(x,t) = \dfrac{1}{2}\int_0^x q(y,t)dy, \quad \pi_{11}(x,t) = \beta_x(x,t)-\beta^2(x,t), \quad \gamma_*(x,t) \text{ by \eqref{eq:Linkage}}.
\]
Deine also
\[ \dfrac{\partial}{\partial x} H_0(x,t) = \dfrac{\partial}{\partial x} \bbmatrix{-\beta(x,t)&-i\dfrac{\pi_{11}(x,t)}{2}\\i\dfrac{\pi_{11}(x,t)}{2} & h_0^{22}(x,t)}
\] 
where  $\dfrac{\partial}{\partial x} h_0^{22}(x,t)$ is real-valued and satisfies the last equation of \eqref{eq:Hi+1HiFs}:
\[ 2i \frac{\partial}{\partial x} h_0^{22}(x,t) = -i \frac{\partial^2}{\partial x^2} \dfrac{\pi_{11}(x,t)}{2} + 2i \beta(x,t) \frac{\partial}{\partial x} \dfrac{\pi_{11}(x,t)}{2}. \]
Then there exists a unique solution $\Phi_*(\lambda,x,t)$, satisfying the following system of equations
\begin{equation} \label{eq:Phi*eqns} \left\{  \begin{array}{llll}
\lambda \sigma_2 \Phi_*(\lambda,x,t) - \sigma_1 \dfrac{\partial}{\partial x}\Phi_*(\lambda,x,t) +
\gamma_*(x,t) \Phi_*(\lambda,x,t) = 0,\\
\dfrac{\partial}{\partial t} \Phi_*(\lambda,x,t) = i \lambda \dfrac{\partial}{\partial x} \Phi_*(\lambda,x,t) + i \dfrac{\partial}{\partial x}[H_0(x,t)] \sigma_1 \Phi_*(\lambda,x,t), \\
\Phi_*(\lambda,0,0) = I.
\end{array} \right. \end{equation}
Notice that for $t=0$ the fundamental matrix $\Phi_*(\lambda,x,0)$ coincides with the fundamental matrix $\Phi_*(\lambda,x)$ considered earlier as the solution of 
\eqref{eq:OutCC}. Evolving this solution using $t$-derivative, we will get the fundamental solution of \eqref{eq:Phi*eqns}. Notice that the system of equations \eqref{eq:Phi*eqns}
is uniquely solvable in $\mathbb R\times[0,T]$. The identity of second mixed partial derivatives 
\[ \dfrac{\partial^2}{\partial x \partial t} \Phi_*(\lambda,x,t) = \dfrac{\partial^2}{\partial t \partial x} \Phi_*(\lambda,x,t)
\]
follows from \eqref{eq:Dgamma*tKdV}, which is identical to \eqref{eq:KdV} due to Corollary \ref{cor:H0toKdV}.
We define the moments, associated with $q(x,t)$ as follows
\[ H_n(x,t) = \int_0^\infty \Phi_*(i\mu,x,t) d\bar\rho(\mu)(i\mu)^n \Phi(i\mu,x,t),
\]
where $d\bar\rho$ is the measure constructed for $q(x,0)$ in the previous section. Then indeed
\[ H_0(x,t) = \int_0^\infty \Phi_*(i\mu,x,t) d\bar\rho(\mu) \Phi(i\mu,x,t).
\]
It follows from the uniqueness of the analytic solution $q(x,t)$ of the KdV equation \eqref{eq:KdV} at the strip $\mathbb R\times[0,T]$. The entry
$h_0^{22}(x,t)$ is chosen from this equality and will satisfy the last equation of \eqref{eq:Hi+1HiFs}, which holds for $\int_0^\infty \Phi_*(i\mu,x,t) d\bar\rho(\mu)\Phi(i\mu,x,t)$
by the construction.

We can define the following operators:
\begin{equation} \label{eq:InvVesselDef}
 \begin{array}{lllll}
B_*(x,t) = \dfrac{1}{2\pi i} \int\limits_\Gamma (\lambda I + A_\zeta)^{-1} B_0 \Phi^*_*(\bar\lambda,x,t) d\lambda, \\
C_*(x,t) = \dfrac{1}{2\pi i} \int\limits_\Gamma \Phi_*(\lambda,x,t) C_0 (\lambda I - A)^{-1} d\lambda, \\
\mathbb X_*(x,t) = I -  \int_\Gamma \int_\Gamma d\lambda(\lambda I + A_\zeta)^{-1} B_0 \dfrac{\Phi_*^*(\bar\lambda,x,t) \sigma_1 \Phi_*(\mu,x,t) - \sigma_1}{\lambda-\mu} C_0 (\mu I-A)^{-1}  d\mu,
\end{array} \end{equation}
where $\Gamma$ is in the sector of regularity of $A, A_\zeta$ as in the previous section.
Along with differential equations \eqref{eq:DB*}, \eqref{eq:DC*}, \eqref{eq:X*Lyapunov} these operators will also satisfy analogues of 
\eqref{eq:DBt}, \eqref{eq:DCt}, \eqref{eq:DXt} as follows:
\begin{eqnarray}
\label{eq:DB*t} \frac{\partial}{\partial t} B_*  & = - i A_\zeta \dfrac{\partial}{\partial x} B_* - i B_* \sigma_1 \dfrac{\partial}{\partial x} H_0,  \\
\label{eq:DC*t} \frac{\partial}{\partial t} C_* u & = i \dfrac{\partial}{\partial x} C_* A u + i \dfrac{\partial}{\partial x}[H_0]\sigma_1 C_* u , \\
\label{eq:DX*t} \frac{\partial}{\partial t} \mathbb X_* v & =  (i A_\zeta B_* \sigma_2 C_* - i B_*\sigma_2 C_* A - i B_* \gamma_* C_* )v ,
\end{eqnarray}
One can argue that these differential equations actually serve as defining ones for the operators $B_*, C_*, \mathbb X_*$ 
with the initial conditions at $t=0$ $\mathbb X_*(x)B(x), C(x)\mathbb X_*(x), \mathbb X_*(x)$ defined in the previous section.
\begin{lemma}[Uniqueness of the moments] \label{lemma:MomentUniquet}Suppose that $H_n(x,t)$ and $\widetilde H_n(x,t)$ are two sequences of moments, which are analytic in $x,t$,
satisfy \eqref{eq:DHntKdV} and $H_n(x,0) = \widetilde H_n(x,0)$ for all $x\in\mathbb R$. Then $H_n(x,t) = \widetilde H_n(x,t)$ for all $x,t$.
\end{lemma}
\Proof It is mmediate from \eqref{eq:DHntKdV}, because for $n=0$, for example,
\begin{multline*}
\dfrac{\partial}{\partial t} H_0(x,0)  = i \dfrac{\partial}{\partial x} H_1(x,0) + i \dfrac{\partial}{\partial x} [H_0](x,0) \sigma_1 H_0(x,0) = \\
= i \dfrac{\partial}{\partial x} \widetilde H_1(x,0) + i \dfrac{\partial}{\partial x} \widetilde [H_0](x,0) \sigma_1 \widetilde H_0(x,0) = 
\dfrac{\partial}{\partial t} \widetilde H_0(x,0) .
\end{multline*}
Continuing in the same manner by induction, we will obtain that 
\[ \dfrac{\partial}{\partial t} H_0^{(n)}(x,0) = \dfrac{\partial}{\partial t} \widetilde H_0^{(n)}(x,0)
\]
and the Lemma follows. \qed
\begin{cor} the following equality holds
\[ C_*(x,t) A^n B(x,t) = H_n(x,t).
\]
\end{cor}
\Proof Since $C_*(x,0) A^n B(x,0) = H_n(x,0)$ and differentiating
\begin{multline*}
 \dfrac{\partial}{\partial t} [C_* A^n B] =
i \dfrac{\partial}{\partial x} C_* A^{n+1} B + + i \dfrac{\partial}{\partial x}[H_0]\sigma_1 C_* A^n B + C^* A^n i A B_x = \\
= i \dfrac{\partial}{\partial x}[ C_* A^{n+1} B]+i \dfrac{\partial}{\partial x}[H_0]\sigma_1 C_* A^n B,
\end{multline*}
we obtain that $C_* A^n B$, $H_n$ satsify the same differential equations with identical initial conditions, so they are equal by Lemma \ref{lemma:MomentUniquet}.
Particularly, for $n=0$ we obtain that $C_*(x,t) B(x,t) = H_0(x,t)$. \qed

From this Lemma it follows that Theorem \ref{thm:vesselmappings} holds, using the same idea of proof, but with the $t$-derivatives
\begin{thm} \label{thm:XX*xt}The following equalities hold:
\begin{eqnarray*}
\mathbb X(x,t) B_*(x,t) = B(x,t),  \quad C_*(x,t)\mathbb X(x,t) = C(x,t), \\
\mathbb X_*(x,t) B(x,t) = B_*(x,t),  \quad C(x,t)\mathbb X_*(x,t) = C_*(x,t).
\end{eqnarray*}
\end{thm}
\Proof We will mimic the proof of Theorem \ref{thm:vesselmappings}.
\[ \begin{array}{llll}
\dfrac{\partial}{\partial t} [C_*\mathbb X] 
& = i \dfrac{\partial}{\partial x} C_* A \mathbb X + i \dfrac{\partial}{\partial x}[H_0]\sigma_1 C_* \mathbb X +
C_* [i A B\sigma_2 C - i B \sigma_2 C A_\zeta + iB\gamma C] \\
& = \text{ using \eqref{eq:DC}, \eqref{eq:Lyapunov}, \eqref{eq:DC*}} \\
& = -i \dfrac{\partial}{\partial x} [C_*\mathbb X] A_\zeta + i \dfrac{\partial}{\partial x}[H_0]\sigma_1 [C_*\mathbb X - C].
\end{array} \]
Plugging here $C$ instead of $C_*\mathbb X$, we will obtain \eqref{eq:DCt} for $C$ and by the uniqueness of the solutions since for $t=0$ it holds that
$C_*(x,0) \mathbb X(x,0) = C(x,0)$ (Theorem \ref{thm:vesselmappings}), we obtain that also $C_*(x,t)\mathbb X(x,t) = C(x,t)$.

In a similar manner one can obtain that
\[ \dfrac{\partial}{\partial t} [\mathbb X_* B] = - i A_\zeta \dfrac{\partial}{\partial x} [\mathbb X_* B] - i B_* \sigma_1 \dfrac{\partial}{\partial x}[H_0].
\]
Since $B_*$ when substituted here instead of $\mathbb X_* B$ satisfies the same equation, and $\mathbb X_*(x,0) B(x,0) = B_*(x,0)$ (Theorem \ref{thm:vesselmappings}),
we obtain that $\mathbb X_*(x,t) B(x,t) = B_*(x,t)$. Then as before
\[ C(x,t)B_*(x,t) = C(x,t)\mathbb X_*(x,t)B(x,t) = C_*(x,t)B(x,t) = H_0(x,t).
\]
The equations $\mathbb X(x,t) B_*(x,t) = B(x,t)$, $\quad C(x,t)\mathbb X_*(x,t) = C_*(x,t)$ follow in the same manner.
\qed

\begin{cor} \label{cor:XxtInverse}The operator $\mathbb X(x,t)$ is invertible for all $x\in\mathbb R$, $t\in[0,T]$ with the inverse $\mathbb X_*(x,t)$.
\end{cor}
\Proof As in Corollary \ref{cor:XxInvertible} we obtain that 
\[ \dfrac{\partial}{\partial x} [\mathbb X(x,t) \mathbb X_*(x,t)] = \dfrac{\partial}{\partial x} [\mathbb X_*(x,t) \mathbb X(x,t)] = 0.
\]
Differentiating with respect to $t$, we obtain that
\[ \begin{array}{lllllllll}
\dfrac{\partial}{\partial t} [\mathbb X(x,t) \mathbb X_*(x,t)] = \text{ \eqref{eq:DXt}, \eqref{eq:DX*t}} \\
= (i A B \sigma_2 C - i B\sigma_2 C A_\zeta+iB\gamma C) \mathbb X_* + 
\mathbb X (i A_\zeta B_* \sigma_2 C_* - i B_*\sigma_2 C_* A - i B_* \gamma_* C_* ) \\
= \text{ Theorem \ref{thm:XX*xt}} \\
= i A B \sigma_2 C_* - i B\sigma_2 C A_\zeta\mathbb X_* + i B\gamma C_*
 + i \mathbb X A_\zeta B_* \sigma_2 C_* - i B \sigma_2 C_* A - i B \gamma_* C_* \\
= \text{ \eqref{eq:Lyapunov}, \eqref{eq:X*Lyapunov}} \\
= i A B \sigma_2 C_* - i B\sigma_2 C [-\mathbb X_* A - B_*\sigma_1 C_*] + i B\gamma C_*
 +i [-A \mathbb X - B\sigma_1 C] B_* \sigma_2 C_* - i B \sigma_2 C_* A - i B \gamma_* C_* \\
= i B\sigma_2 C B_*\sigma_1 C_* + i B\gamma C_*
 - i B\sigma_1 C B_* \sigma_2 C_* -  i B \gamma_* C_* \\
= i B[\sigma_2 H_0 \sigma_1 + \gamma - \sigma_1 H_0 \sigma_2 -\gamma_* ] C_* = \text{\eqref{eq:Linkage}} \\
= 0.
\end{array} \]
Similarly, one shows that $\dfrac{\partial}{\partial t} [\mathbb X_*(x,t) \mathbb X(x,t)] = 0$. Thus the operator $\mathbb X(x,t)\mathbb X_*(x,t)$
is analytic and has zero derivatives with respect to $t$, and $x$. Thus it is constant. Since at $x=t=0$ it is identity, the result follows. \qed
\begin{thm} \label{thm:XinvForq}
Suppose that $q(x,t)$ is a solution of \eqref{eq:KdV} on $\mathrm R\times[0,T]$, then there exists a KdV vessel, realizing $q(x,t)$ on
$\Omega$, so that $\mathrm R\times[0,T]\subseteq\Omega$.
\end{thm}
\Proof For $q(x,0)$ we construct the SL vessel $\mathfrak{V}$. Since the solution $q(x,t)$ exists on $\mathrm R\times[0,T]$ the fundamental matrix 
$\Phi_*(\lambda,x,t)$ of \eqref{eq:OutCC} exists for all $(x,t)\in\mathrm R\times[0,T]$. We can also define the ingredients of the ``inverse vessel''
$B_*(x,t), \mathbb X_*(x,t), C_*(x,t)$ by formulas \eqref{eq:InvVesselDef}. Then by Corollary \ref{cor:XxtInverse} $\mathbb X_*(x,t)$ is the inverse of $\mathbb X(x,t)$ on 
$\mathbb R\times[0,T)$. So, $\Omega$ includes the set $\mathrm R\times[0,T]$. \qed

\section{\label{sec:Remarks}Remarks}
\noindent\textbf{1.} A next step, related to this research, would be a development of a similar theory for locally integrable functions $q(x)$ on $\mathbb R$. Using approximations
by analytic functions similar results should be obtainable.

\noindent\textbf{2.} In the case $q(x)$ is analytic and satisfies $\int\limits_{-\infty}^\infty (1+|x|) q(x) dx<\infty$ \cite{bib:FaddeyevII} one obtains that the fundamental matrices
$\Phi_*(i\mu,x), \Phi_*(-A_\zeta,x)$ are UNIFORMLY bounded, and as a result the operator $\mathbb Y^*(x)$ is uniformly bounded as well. 
Thus there exists $T = T_x$ - the same for all $x$, so that the operator $\mathbb X(x,t)$ is invertible. It follows that there exists
a local solution on $\mathbb R\times[-T,T]$ of the KdV equation  \eqref{eq:KdV}, which is a very well known result for the KdV equation \cite{bib:GGKM, bib:FaddeyevII}.

\noindent\textbf{3.} These vessel constructions are intimately related to the theory of systems \cite{bib:Staffans}. Their relations and interplay is left for a future work. 

\bibliographystyle{alpha}
\bibliography{../../biblio}

\end{document}